\documentclass[11pt]{article}
\usepackage{amsfonts,amsmath,amssymb,amsthm,amscd,epsfig,epic,eepic}
\usepackage{graphics}

 \evensidemargin=\oddsidemargin

\makeatletter
\gdef\thmhead@plain#1#2#3{%
  \thmname{#1}\thmnumber{\@ifnotempty{#1}{ }#2}%
  \thmnote{ {\mdseries#3}}}
\let\thmhead\thmhead@plain
\makeatother
\theoremstyle{plain}
\newtheorem{theorem}{Theorem}[section]
\newtheorem{lem}[theorem]{Lemma}
\newtheorem{proposition}[theorem]{Proposition}
\newtheorem{corollary}[theorem]{Corollary}

\newtheorem*{lemi}{Lemma}

\newtheorem{claima}{Claim}
\newtheorem{claimb}{Claim}
\newtheorem{claim}{Claim}
\newtheorem{pro}{Proposition}
\newtheorem*{proi}{Proposition}

\newtheorem{assumption}{Assumption}
\newtheorem{thm}{Theorem}

\newtheorem*{crr1}{Corollary 1}

\newtheorem*{conj}{Conjecture}
\newtheorem*{proprietes}{Properties}

\theoremstyle{definition}
\newtheorem{definition}[theorem]{Definition}
\newtheorem{remark}[theorem]{Remark}
\newtheorem{notation}[theorem]{Notation}

\theoremstyle{remark}

\def\alinea#1{\hfill\break%
  \hbox to \parindent{\hss{\upshape{\bf #1)}}\enspace}\ignorespaces}
\def\bul{\hfill\break\hbox to\parindent{\hss$\bullet$\enspace}\ignorespaces}

\def\from{\colon}

\def\from{\colon}
\def\mod{\operatorname{mod}}

\def\wc{\omega_{comb}}

\def\H{H}
\def\C{\mathbf{C}}
\def\D{\mathbf{D}}

\def\N{\mathbf{N}}

\def\R{\mathbf{R}}
\def\S{\mathbf{S}}

\def\AA{\mathcal{A}}\def\BB{\mathcal{B}}
\def\CC{\mathcal{C}}

\def\DD{\mathcal{D}}

\def\LL{\mathcal{L}}

\def\PP{\mathcal{P}}
\def\PPc{\PP_{\omega Crit(c)}}
\def\PPo{\PP_{\omega Crit(c_0)}}

\def\RR{\mathcal{R}}

\def\Cap_#1{\bigcap\limits_{#1}}
\def\Cup_#1{\bigcup\limits_{#1}}
\def\ol{\overline}

\def\wt{\widetilde}

\def\eps{\epsilon}
\def\Sun{S^1}
\newcommand\TTT{{\cal T}}
\newcommand{\veps}{\ensuremath{{\underline{\epsilon}}}}

\makeatletter
\def\cqfdsymb{\relax\protect\ifmmode\else\unskip\nobreak\fi
\quad\hfill$\bgroup
\vcenter{\hrule\hbox{\vrule\@height.6em\kern.6em\vrule}\hrule}\egroup$}
\def\cqfd{\cqfdsymb  \endtrivlist}
\gdef\rom#1{\leavevmode\skip@\lastskip\unskip\/%
        \ifdim\skip@=\z@\else\hskip\skip@\fi{\normalshape#1}}
\makeatother
\begin{document}
\title{{\bf Bounded critical Fatou components are  Jordan domains, for polynomials.}}

\author{{{\sc P. Roesch} , {\sc Y.Yin}\ }\thanks{Research partially supported by the ANR ABC and the NSF of China}
\\
{\small Institut de Math\'ematiques de Toulouse, France\,; Fudan University, China}}
\date{\today}
\maketitle

\begin{abstract} We prove that the boundary of  any bounded Fatou component 
for  a polynomial is a  Jordan curve,  except maybe for  Siegel disks.\end{abstract}

\section{Introduction}
The Riemann sphere, viewed as  a dynamical space on which act the  rational maps, is divided into two sets  : 
 the Fatou set  $F(f)$ and the Julia set $J(f)$. The dynamics on the Fatou  set is well understood but is chaotic on the Julia set. In this article we will restrict ourself to the action of polynomials. The connected component of the Fatou set
   containing $\infty$, usually called $B(\infty)$,  is distinguished\,:  its boundary is exactly  $J(f)$.
It can have a very 
complicated  topology  (for instance it can contain continua of the form of a  "hedgehog" discovered by Perez-Marco).
 However, the  bounded Fatou components have  very simple  boundary : 
\begin{thm}\label{th:fatouboundaryJordan} 
If $f$ is a polynomial, any  bounded Fatou component, which is not a Siegel disk, is 
a Jordan domains {\rm(}{\it i.e.}
a disk with Jordan curve boundary{\rm)}.
\end{thm}
Concerning  {\it  Siegel disks}, {\it i.e.}  Fatou components  on which the map is conjugate to  irrational rotations,  the question is still open. One conjectures :
\begin{conj} For a polynomial, every bounded Fatou component is a Jordan domain.
\end{conj}
This conjecture is supported by the work of~\cite{PZ} and  Shishikura announced it for high type Siegel disks in degree $2$.
Note that Ch\'eritat recently constructed a  holomorphic map  (that is not a polynomial) defined in an open disk  of $ \C$ containing  the closure of a fixed Siegel disk,   whose boundary is a pseudo circle, in particular it is not locally connected (see~\cite{Ch}).

\vskip 1em 
As a consequence of the proof of Theorem~\ref{th:fatouboundaryJordan} we obtain a ``description from inside'' of the dynamics on  the component containing $U$ of the filled Julia set, $ K(f)=\C\setminus B(\infty)$\,:

\begin{thm}\label{th:description}
Let $f$ be a polynomial and $U$  a periodic bounded Fatou component which is not a Siegel disk. Every point of $\partial U$ is the landing point of at least one external ray. Moreover if  $K_U$ denotes the connected component  containing $U$ 
of  the filled Julia set $K(f)$, then  $\displaystyle K_U=\overline U\cup \bigcup_{t\in \S^1} L_t$ 
where the sets $L_t$ are  ``limbs'' sprouting out of $\overline U$ with  the following properties\,:  
 \begin{enumerate}
\item $L_t$ is connected\,;
\item $L_t$ intersects $\overline U$ at exactly one point called $\gamma_U(t)$\,;
\item $L_t\neq \gamma_U(t)$  if and only if $L_t$ either contains a critical point or is eventually mapped to a limb $L_{t'}$ containing a critical point.
\end{enumerate}
\end{thm}

\begin{crr1}\label{c:lcfails}If $J(f)$ is connected, the only point of $\partial U$, where the local connectivity  of $J(f)$ can fails  are the  eventually periodic ones,   where a comb might be attached by one point to the boundary. 
\end{crr1}
\vskip 1 em \noindent{\it Acknowledgment}: The first author would like to thank Tan Lei for many helpful comments.

\subsection{Overview}
Theorem~\ref{th:description} follows from the proof of Theorem~\ref{th:fatouboundaryJordan} which is a consequence of the following.

\begin{thm}\label{th:fatoulc} Let $f$ be a polynomial and $U$ a  bounded 
Fatou component. Assume that  $J(f)$ is connected and that $f$ fixes $U$  with  exactly  one critical point in $U$, then  $\partial U$ is locally connected.
\end{thm}
\begin{lem}\label{l:thm3}Theorem~\ref{th:fatoulc} implies Theorem~\ref{th:fatouboundaryJordan}.
\end{lem}
\proof
Let $f$ be a polynomial and $U$ a bounded Fatou component which is not a Siegel disk.

Denote by $K_U$ the connected component of $K(f)$
 containing $U$.  There exist  neighborhoods $X,X'$ of $K_U$, an integer $r$ and a polynomial $g$ with connected Julia set, such that  the iterate $f^r$ of $f$maps $X$ to $X'$ and $f^r : X\to X'$ is conjugated  by a quasi-conformal homeomorphism to $g$ on a disk $D(0,R)$  with large $R$. This follows from the classical theory of polynomial-like mappings applied to some disk of low potential (for the Green function) containing $K_U$.
Let $V$ be the image of $U$ by the conjugacy\,; it is equivalent to prove that $\partial U$ or $\partial V$ is a Jordan curve. 

By Sullivan's non-wandering Theorem, $V$ is mapped by $g$ to some Fatou component $W$ which is fixed by some iterate $h=g^k$. Moreover, this Fatou component $W$ contains a critical point so cannot be a Siegel  disk\,: a disk where the dynamics is conjugated to an irrational rotation. Hence by the classification result of the Fatou components, there is a fixed point $p=h(p)$ in $\ol W$ which is either attracting when $p\in W$,  or parabolic when  $p\in\partial W$.  By a surgery procedure $h$ is quasi-conformally conjugated on a neighborhood of the Julia set to a polynomial having only one critical point in the  Fatou component $Y$  corresponding to $V$. 
This surgery is done  for instance in  Theorem~5.1 of~\cite{CG} for  the attracting case and, in the  parabolic case in   Proposition~6.8 of~\cite{McMautom} (where  general case is in fat done).

Thus,  since the Julia set of an iterate (here  $g^k$) is equal to the original Julia set (here of $g$), 
we have  just proved that there is  an homeomorphism between  neighborhoods of the boundary  of the corresponding Fatou components,
$U$ and   $Y$. Therefore, applying Theorem~\ref{th:fatoulc} to the polynomial $h$ and the Fatou component $Y$, implies  the local connectivity of $\partial U$.

Then it is classical to prove that a it is a Jordan domain using the maximum principal.
The Riemann map $\Phi $, 
from $\D$ to $U$ (which is a disk), extends continuously to $\ol \D$ since  $\partial U$ is locally connected (by  Carath\'eodory's Theorem). Therefore, $\partial U$ is the   curve $\Phi(\S^1)$. If it is not simple, there exists 
$t,t'$ such that $\Phi(e^{2i\pi t})=\Phi(e^{2i\pi t'})$ so that the curve $C=\Phi([0,1]e^{2i\pi t})\cup  \Phi([0,1]e^{2i\pi t'})$ is in $\ol U$ and bounds points that are attracted by  $\infty$,  which contradicts   the maximum principle.
 \cqfd

\subsection{About the Proof of Theorem~\ref{th:fatoulc}}

We find connected neighborhoods for  the  points  $x\in \partial U$  as pieces of the complement of some backward iterated graph (classically called {\it puzzle pieces}). These graphs are  the union of internal and external rays as well as equipotentials  (in the parabolic case we use "parabolic rays"). To prove that the diameter of the pieces  tends to zero we use different techniques depending on the point $x$ of $\partial U$\,:
\begin{itemize}
\item If $x$ is eventually periodic (see section~\ref{s:kiwi}) it might be the attach-point of some other part 
of the Julia set. We consider the dynamics on the intersection of the closure 
of the puzzle pieces  containing this periodic point (following the ideas of  Kiwi in~\cite{Kiwi}) and prove that two external rays, landing  at $x$, separate
this continuum from  the closure of $U$. This is the only case where the intersection of the puzzle pieces might not be reduced to one point and then    we might have a renormalizable map.
\item
If $x$ {\it combinatorially accumulates } an  eventually periodic  point $y$, meaning   for the topology generated by the set of the puzzle pieces,  we consider the first entrance time in the periodic  nest\,: the first time the orbit enters in each puzzle piece of the nest around $y$ (see section~\ref{s:techniques}). 

If $y$ is repelling, it is easy to find a non degenerate annulus between consecutive puzzle pieces and to pull it back 
with bounded degree around $x$.  If $y$ is parabolic,  we use some distortion properties on the enlarged puzzle pieces.
\item
The last case to consider is when  $x$ does not combinatorially accumulate eventually periodic points. Here, the recurrence of the critical points plays a fundamental role.
We consider  the  combinatorial accumulation of the different critical points on themselves\,:
\begin{itemize}
\item The ``non-recurrent case''  (see section~\ref{s:techniques}) corresponds to the situation where there is at least a critical point  (in the combinatorial accumulation of $x$) whose combinatorial accumulation does not contain critical points. This case works as before, excepted that for finding a non degenerate annulus, where we need to look at deeper puzzle pieces in the nest.
\item
In the  ``recurrent case''  we have to distinguish between two strength of recurrences\,: the ``reluctantly recurrent'' case is dealt in section~\ref{s:techniques} using long iterates of bounded degree, and the ``persistently recurrent'' case 
(see section~\ref{s:persrec}) where  we have to introduce the {\it enhance nest} introduced by  Koslovski, Shen and vanStrien in~\cite{KSS}. It is a  double sub-nest $(K_n,K'_n)$ of the nest of a critical point and that has the property that $K'_n\setminus K_n$ avoids the orbit of the critical points, that the time to go from  $K_n$ to $K_{n-1}$ is more  than half of the time to go from $K_n$ back to $K_0$,  whereas the degree of the map from $K_n$ to $K_{n-1}$ is uniformly  bounded.  To prove that  the diameter of the puzzle of this nest  tends to zero, we 
use  the Covering Lemma of Kahn and Lyubich (see~\cite{KL}).
The arguments are similar to those provided in  the proof of  the Branner-Hubbard conjecture  (see~\cite{QY} and also~\cite{KS}), excepted that we do not have non degenerate annuli here.

To have a different point of view on the enhanced nest we recommend highly to read~\cite{TY} (in the  case of a unique critical point) and~\cite{TLP} (in the case of several critical points).
\end{itemize}\end{itemize}

Excepted for section~\ref{s:KU}, the rest of the article is devoted to  the proof of  Theorem~\ref{th:fatoulc}. 
\section{Puzzles and parabolic techniques}
\subsection{Notations}
We fix a    polynomial  $f$  of degree $D$ and assume that it has a  connected filled-in Julia set $K(f)$. We consider  a   bounded Fatou component $U$, that is fixed by $f$ and  which contains exactly one critical point.

Since $K(f)$ is connected,  its complement $B(\infty)\cup\{\infty\}$ is simply connected. Denote by $\Phi_\infty:\ol\C\setminus \ol\D\to B(\infty)\cup\{\infty\}$ the Riemann map such that
$\Phi_\infty(\infty)=\infty$ and which is tangent to  the identity near $\infty$.
It conjugates $f$ to $z\mapsto z^D$. The bounded Fatou component $U$ is a topological disk, let $c_0$ be the unique  critical point of $f$ in $U$ and denote by $\Phi_U:\D\to U$ a Riemann map such that $\Phi_U(0)=c_0$. 

In {\it the attracting case},  $c_0$ is fixed (it is the super-attracting fixed point of $U$).
We choose the map   $\Phi_U$ such  that  $\Phi_U'(0)=1$. Therefore,  
 $\Phi_U$ conjugates $f$ to  $z\mapsto z^d$ where $d$ is the degree of the critical point $c_0$.

In the {\it parabolic case}, there is a fixed point $p$ on $\partial U$ and every point of $U$ tends to $p$ under the iteration of $f$. 
We choose  $\Phi_U$ such  that $\Phi_U(v_d)=f(c_0)$ where   $v_d=\frac{d-1}{d+1}$.
Therefore,  
$\Phi_U$ conjugates $f$ to  the  Blaschke model $B(z)=\frac{z^d+v_d}{1+v_dz^d}$ on $\D$ and maps $p$  to $1$.

\subsection{Rays and equipotentials} 
\underline{\it In the attracting case}\,:
\vskip 0.5em
\begin{definition}The  {\it external }, resp. {\it  attracting internal,  ray}  of angle $\theta$  is the set 
$$R_\infty(\theta)=\Phi_\infty(\{e^{2i\pi \theta+t}\mid t>0\}), \hbox{ resp. } R_U(\theta)=\Phi_U(\{e^{2i\pi \theta+t}\mid t<0\}).$$
The  {\it external}, resp. {\it internal,  equipotential} (in the attracting case), of potential  $v>0$  is the set 
$$E_\infty(v)=\Phi_\infty(\{e^{2i\pi \theta+v}\mid \theta\in\R\}), \hbox{ resp. } E_U(v)=\Phi_U(\{e^{2i\pi \theta-v}\mid \theta\in\R\}).$$
\end{definition}

\begin{pro}[(Douady, Hubbard, Sullivan, Yoccoz)]\label{th:extlanding}
If $t$ is rational, the internal ray $R_U(t)$ lands at a (eventually) periodic point which is either repelling or parabolic. Moreover, there is at least one  external ray landing at this point. All these rays (internal and external)  have the same rotation number.  
\end{pro}
We say that a $q$ cycle for $f^k$  of rays $R_0, \ldots R_{q-1}$
landing on a common $k$ periodic point $z$ and numbered in the counter clockwise order
around $z$ defines \emph{combinatorial rotation number} $p/q$, $(p,q)=1$
iff $f^k(R_j) = R_{(j+p)\mod q}$.

\vskip 1em\noindent\underline{\it In the parabolic  case}\,:
\vskip 0.5em 
We recall a definition of parabolic rays given in~\cite{PR}. 
We  first construct parabolic rays in the disk for the model map $B$ and  then we pull them back by the conjugacy.

Let $z_\emptyset=0$ and $T_\emptyset := B^{-1}([0,v]) = \cup_{j=0}^{d-1} [0,z_j]$, 
where $z_j={v_d}^{1/d}e^{i\pi/d}\omega^j$ and $\omega=e^{2i\pi /d}$.
Moreover let $T_j$ denote the connected component of
$B^{-1}(T_\emptyset)$ containing $z_j$.  
Define recursively on $n\in\N^*$ and for each
$(\eps_1,\eps_2,\ldots,\eps_n)\in{\{0,1,\ldots,(d-1)\}}^n$ the
point $z_{\eps_1,\eps_2,\ldots,\eps_n}$ as the unique point of the
preimage $B^{-1}(z_{\eps_2,\ldots,\eps_n})$ belonging to
$T_{\eps_1,\eps_2,\ldots,\eps_{n-1}}$. And define
$T_{\eps_1,\eps_2,\ldots,\eps_n}$ to be the connected component of
the preimage $B^{-1}(T_{\eps_2,\ldots,\eps_n})$ containing
$z_{\eps_1,\eps_2,\ldots,\eps_n}$.

Define for each $n$ a ($d$-adic) portion of  tree
$\TTT_n:=\cup_{k=0}^nB^{-k}(T_\emptyset)$, so that
$$\TTT_n = \TTT_{n-1}\cup\bigcup_{(\eps_1,\eps_2,\ldots,\eps_n) \\
\in{\{0,1,\ldots,(d-1)\}}^n} T_{\eps_1,\eps_2,\ldots,\eps_n}.
$$

Moreover define an infinite $d$-adic tree $\TTT:=\cup_{k=0}^\infty
B^{-k}(T_\emptyset)$ with boundary (in) $\S^1$. 
\begin{definition}[{\rm (}Parabolic rays for $B${\rm)}] For $\veps\in\Sigma_d$,  where $\Sigma_d:={\{0,1,\ldots,d-1\}}^{\N}$, define the {\it parabolic ray }$R_\veps$ as the
minimal connected subset of $\TTT$ containing the sequence of
points $(z_{\eps_1,\eps_2,\ldots,\eps_n})_{n\in \N}$ (enterpreting $n=0$ as $z_\emptyset$).
\end{definition}
For each $0\leq j< d$, let  $S_j$ be the open sector  spanned by the interval
$I_j=[\omega^j,\omega^{j+1}]\subset \S^1$, i.e. the interior of the convex hull of the union of $I_j$
and $0$. The map $B$  is univalent from $S_j$ onto $\D\setminus[v,1]\supset
T_\emptyset$. Its boundary arcs $[0,\omega_j]$ and
$[0,\omega_{j+1}]$ are each mapped (homeomorphically) onto  the
forward invariant arc $[v,1]$. As $z_j\in S_j$ for each $j$ it
easily follows by induction, that $\TTT_n\cap S_j$ is connected
and for any $n$ and any $(\eps_2,\ldots,\eps_n)$ contains both the
point $z_{j,\eps_2,\ldots,\eps_n}$ and the set
$T_{j,\eps_2,\ldots,\eps_n}$. Consequently $S_j\cup\{0\}$ contains
also any of the rays $R_\veps$ with $\epsilon_1=j$.

Consider the attracting Fatou coordinate $ \Phi_+\from \D \to \C$ for
$B$ on $\D$, normalized by $\Phi_+(0) = 0$. For any $\veps\in \Sigma_d$, $\Phi_+$
 maps $R_\veps$
homeomorphically onto $\R_-=]-\infty,0]$, and has a degree $d$
critical point at the preimage $z_{\eps_1,\eps_2,\ldots,\eps_n}$
of $-n$ for each $n\geq 0$.  The  extended ray  $\widehat R_\veps:=R_\veps\cup [0,1[$ is mapped homeomorphically to $\R$ by
$\Phi_+$. Denote by$\widehat R_\veps (t)$ the point
$\Phi_+^{-1}(t)\cap \widehat R_\veps$ and define similarly
$R_\veps(t)$ for $t\leq 0$. By construction,
$$ \forall\;\veps\in\Sigma_d :
B(\widehat R_\veps(t))=\widehat R_{\sigma(\veps)}(t+1), \; \hbox{ where } \sigma(\epsilon_1,\epsilon_2,\ldots) =
(\epsilon_2,\epsilon_3,\ldots)  \; \hbox{ is the shift map on } \Sigma_d. $$
Note  that  every ray $R_\veps(t)$ lands at some point $z_\veps\in \S^1$, as $t\to -\infty$.

\begin{definition}[{\rm (}Parabolic rays in $U${\rm)}] For $\veps\in\Sigma_d$, 
define the parabolic ray of argument $\veps$ in $U$
as $R_U[\veps]:=\Phi_U^{-1}(R_\veps)$ and the extended  ray as
$\widehat R_U(\veps):= \Phi_U^{-1}(\widehat R_\veps)$.
\end{definition}
\noindent
By construction,
$f(\widehat R_U[\veps](t))=\widehat R_U[\sigma(\veps)](t+1),  \forall\;\veps\in\Sigma_d 
$.
\vskip 1em
The correspondence  given by the formula $ \theta = \sum_{n=1}^\infty \frac{\epsilon_n}{d^n}$, between angles and itineraries ($\veps$), is a bijection for non d-adic angles.
Therefore, we will denote by    $R_U(\theta)$ the  ray of angle $\theta$, and by 
$\widehat R_U(\theta)$ the extended ray, as soon as $\theta$ is not d-adic.
We say that  the parabolic ray  $\widehat R_U[\veps]$ converges if $\widehat R_U[\veps](t)$ admits a limit when $t$ tends to $-\infty$.

\begin{thm}\label{simplelanding}
For any (pre-)periodic argument $\veps\in\Sigma_d$,
i.e.~$\sigma^k(\sigma^l(\veps))=\sigma^l(\veps)$,
the parabolic ray $R=R_U(\veps)$ converges to an $f$ (pre-)periodic point
$z\in\partial U$ with $f^k(f^l(z))=f^l(z)$.
Furthermore, if $\veps$ is periodic (i.e.~$l=0$), let $k'$ denote the exact period of $z$ and let $q=k/k'$.
Then the ray $R$ defines a combinatorial rotation number $p/q$, $(p,q)=1$ for $z$.
The periodic point $z$ is repelling or parabolic with multiplier $e^{i2\pi p/q}$.
Moreover any other ray in $U$ landing at $z$ is also $k$ periodic and
defines the same rotation number.
\end{thm}
This is a standard result which in its initial form is due to
Sullivan, Douady and Hubbard, for external rays of polynomials,  see also \cite[Th. A and Prop. 2.1]{Petersen1}.

\begin{definition} Consider for $v>0$, the closure of the union of the curves   $\Phi_U(\Phi_+^{-1}(\{ z \mid \Re e(z)=\ln v/\ln d\}))$. Denote by $E_U(v)$ the component containing $p$.  For $0<v\le 1$,  there is only one  component  in this set, however for  $v> 1$, this set contains several components.
\end{definition}
For $v \in \R^{+*}\setminus \{1/d^k\mid k\ge 0\}$,  $E_U(v)$  is a simple closed  curve. Moreover,   it satisfies the same relation as in the attracting case, namely that $f(E_U(v))= E_U(dv)$ for $v>0$.

\begin{lem}\label{l:aboutdistray}If $\theta\neq\theta'$ the rays $R_U(\theta)$ and $R_U(\theta')$  do not land at the same point (if they land).
\end{lem}
\proof Assume that the two rays land at the same point, then the union of their closure defines  a closed curve $\gamma$. There are points of the Julia set in the bounded connected component of $\C\setminus \gamma$\,: take for instance the landing point of dyadic angles that  are between $\theta$ and $\theta'$. This contradicts the maximum principle.\cqfd 

\subsection{Puzzles}

Let us define the puzzle associated to a graph defined by a periodic angle $\theta$. Then taking the appropriated preimage will  lead to the  definition for a pre-periodic angle.
If $\theta$ is a $k$-periodic angle by multiplication by $d$, but not fixed, the internal ray  
$R_U(\theta)$ is well defined (in the case of parabolic ray) and lands at  a periodic point
 $z(\theta)$ of $\partial U$.
Denote by $R_\infty(\xi)$ an external ray landing at $z(\theta)$ (Proposition~\ref{th:extlanding}). It as the same period 
and same rotation number as $R_U(\theta)$.

\begin{definition}\label{d:graph} Let $\Gamma(\theta)$ be the graph $$ \Gamma(\theta)= \left(\bigcup_{j\ge 0} R_U(d^j \theta)\cup z(d^j \theta)\cup R_\infty(D^j \xi)\right)\cup E_\infty(1)\cup E_U\left(\frac{1}{d^{k-1}-1}\right).$$
\end{definition}
We have the following stability property\,:  $f(\Gamma(\theta))\cap \Omega \subset 
\Gamma(\theta)$, where  $\Omega$ is the connected component of $\ol \C \setminus  E_\infty(1)\cup E_U(\frac{1}{d^{k-1}-1})$ containing $\partial U$ (or a portion  of it).

The definition of $\Gamma(\theta)$ (namely the value of the equipotential in $U$) is chosen so that two rays of $\Gamma(\theta)$ in $U$ do not touch in $\Omega$.
\begin{definition} Let $\Gamma_n(\theta)$ be the graph $ f^{-n}(\Gamma(\theta))$. The connected components of $\ol \C \setminus \Gamma_n(\theta)$ intersecting $\partial U$
 are called {\it puzzle pieces of depth $n$}.
For $y\in \partial U$, let us  denote by $P_n(y)$ the puzzle piece containing $y$, when it exists. We extend this notation to the parabolic point $p$, denoting by $P_n(p)$ the puzzle piece containing $p$ in its closure.
\end{definition}
The next properties follow from the definition of the puzzle pieces. The proof is classical (see~\cite{Ro1} for instance).
\begin{lem}
\begin{itemize}\item
For every point $y\in\partial U$,  which is not eventually mapped to $p$, there exists a strictly periodic angle $\theta$ such that 
 the orbit of $y$ does not intersect  $\Gamma(\theta)$. Hence  $P_n(y)$ is well defined for every $n\ge 0$. Moreover, we can chose $\theta$ such that the orbit of the  critical points  also never cross  the graph $\Gamma(\theta)$ except maybe at the pre-images of
 $p$\,;
\item
The puzzles pieces are either nested or disjoint\,;
\item The image $f(P_n(y))$ is the puzzle piece $P_{n-1}(f(y))$\,;
\item If $y\in \partial U$, the boundary of $P_n(y)$ intersects $U$ along exactly two portions of internal rays and a portion of equipotential that might be touching $\partial U$ at a pre-image of $p$\,;
\item The boundary of two nested distinct puzzle pieces can only touch at  a point of $\partial U$ which is   fixed  or eventually periodic.\end{itemize}\end{lem}

\begin{notation}\label{n:critend} Let $End(y)$ denote the set $\bigcap_{n\in \N} \ol P_n(y)$.  Let $Crit$ denotes the union of the critical points of $f$ outside of $\bigcup_{n\in \N}f^{-n}(U)$  and call {\it critical  ends} the  ends  containing at least a critical point, {\it i. e.} the ends $End(c)$, for $c\in Crit$.
\end{notation}

\begin{assumption}\label{a:crit} Let $Y$ be a finite set of points. We can always assume that there   is no critical point in $P_0(y)\setminus End(y)$ for every point $y\in Y$. 
\end{assumption}
Indeed, one can  replace  $\theta$ by one of its iterated pre-image for instance. We will use this mainly with $Y=Crit\cup \{x\}$ or some iterate of it where $x$ is the point we focus on.
\begin{remark}\label{r:enddegree} Under this assumption we can define the degree on a  critical end $End(c)$ as the degree on $P_0(c)$. It does not depend on the choice of the critical point in the end. For any $y$, we denote by $\delta(y)$ the degree of $f$ on $P_0(y)$.
\end{remark}
\begin{lem}\label{l:connexite} Let $V$ be the closure of a puzzle piece. Then $V\cap\partial U$ is connected.
\end{lem}
\proof  One easily sees by induction that the boundary of $V$  intersects  $U$  along the closure of two rays $R_U(t_1), R_U(t_2)$.  Therefore, $V\cap\partial U$ is the intersection of the sets   $S_n =\ol{ \Phi_U( \Delta_n)}$  where $\Delta_n =\{z\in \C\mid \vert z\vert \in[e^{1/n},1[,\ \arg(z)\in [2\pi t_1,2\pi t_2] \}$. Hence $V\cap\partial U$ is a connected set as the decreasing  intersection of the compact connected sets $S_n$. 
 \cqfd
 
 \begin{corollary} Let $I_n(x)=\ol{P_n(x)}\cap \partial U$. It suffices to prove that $\bigcap_{n\in\N} I_n(x)=\{x\}$  in order to show that  $\partial U$ is locally connected at  the point $x$.
 \end{corollary}
 
\begin{definition}\label{d:changementdegraphe} If $x$ lies on some graph $\Gamma_n$, one can define $ I_n(x)$ by ${Q_n(x)}\cap \partial U$ where $Q_n(x)$ is the closure of the union of the pieces containing $x$ in their closure. 
\end{definition}
\begin{lem}\label{l:changementdegraphe}Let $x\in \partial U$, for any graph $\Gamma(\theta)$ of the Definition~\ref{d:graph} the set $I_n(x)$ is well defined. Moreover, the property $\bigcap_{n\in\N} I_n(x)=\{x\}$ does not depend on the graph.
 \end{lem}
 \proof Consider two graphs $\Gamma, \Gamma'$\,; denote  $P_n, P_n'$ the pieces of the   puzzles constructed from these graphs.  For any $n\ge 0$, there exists $k(n)$ such that $P'_{k(n)}(x)$ (resp. $Q'_{k(n)}(x)$) is included in 
$\ol{P_n(x)}$ or in  $\ol{Q_n(x)}$. Indeed, by construction  the intersection of $P'_k(x)$ with $U$ 
is the sector  delimited  by two rays  $R_U(t_k), R_U(s_k)$\,; the difference $s_k-t_k$ tending to $0$.
Therefore  $P'_k(x)$  will be included in $P_n(x)$   for large values of $k$  (or in $Q_n(x)$  if $x$ lies on $\Gamma'$ or on some inverse image of $\Gamma'$). 
 \cqfd
\section{The ``periodic'' case}\label{s:kiwi}
The Julia set $J(f)$ might contain  {\it Cremer points, i.e.} periodic points with multiplier $e^{2i\pi t}$ for some irrational $t$. Recall that a  Cremer point can not be on $\partial U$ by  a result of Goldberg and Milnor (see~\cite{GM}). But there is some continuum,  containing the Cremer point,   that  possibly  touches $\partial U$ and  is then  the cause of  non local-connectivity at that point.

Since  our  graphs  $\Gamma_n(\theta)$ will  not cut such a continuum (it  is indecomposable), the corresponding sequence of puzzle pieces does not shrink to a point, {\it i.e.} $End(x)\neq\{x\}$.
However we prove that  $End(x)$ is attached by (at most) one point to $\partial U$, which implies that $\cap_n I_n(x)=\{x\}$\,: 

\begin{proposition}\label{p:kiwi} Let $x \in \partial U$ be  a fixed point of  $f$ and let $\Gamma(\theta)$ be a graph (as in Definition~\ref{d:graph}) with  $\theta$  any periodic angle that is not fixed by multiplication by $d$. We can define $P_n(x)$ as the unique puzzle piece containing $x$ in its closure and $End(x):=\cap\ol { P_n(x)}$ . Then
 \begin{itemize}
\item  either $End(x)=\{x\}$\,;
\item  or there exist two  rays $R_\infty(\zeta),R_\infty(\zeta')$ which land at $x$ and separate  $\overline U$ from $ {End(x)\setminus \{x\}}$.\end{itemize}
\end{proposition}

 \begin{corollary}\label{c:kiwi}For any eventually periodic point  $x\in \partial U$, $\cap I_n(x)=\{x\}$. In particular, $\partial U$ is locally connected at $x$. \end{corollary}
\proof 
 Let $x\in \partial U$ be  a fixed point of $f$, take some graph $\Gamma$ as defined in the Proposition. Suppose that $End(x)\neq \{x\}$, consider the curve $\gamma=\ol R_\infty(\zeta)\cup\ol R_\infty(\zeta')$. 
Let $V_n$ be  the closure of the connected component of $P_n(x)\setminus \gamma$ intersecting with $U$.  Note that $\cap _n(V_n\cap \partial U)= I(x)$. Since $\gamma$ separates $\overline U\setminus \{x\}$ from $ End(x)=\cap \ol P_n(x)$, it follows that $\bigcap V_n=\{x\}$. So that $\partial U$ is locally connected at $x$ since $(V_n\cap \partial U)$ is a decreasing sequence of connected neighborhoods of $x$ in $\partial U$.

If $x$ a periodic point on $\partial U$, we consider an iterate $g$ 
of $f$ such that $x$ is a fixed by $g$. Then applying previous result to $g$, we obtain  a sequence of connected neighborhoods of $x$ in $\partial U$ (the Julia set and the basin are the same for $f$ and the iterate $g$ of $f$).

If $x$ is an eventually periodic point on $\partial U$, we pullback by some iterate of   $f$ the previous sequence 
of connected neighborhoods of the periodic point in the orbit of $x$.
\cqfd

\vskip 1em \noindent{\bf Proof of the Proposition~\ref{p:kiwi} :}
Let $\theta$ be any periodic angle that is not fixed by multiplication by $d$.
Consider the sequence of puzzle pieces  $P_n(x)$ defined by the graph $\Gamma(\theta)$ and containing $x$ in their closure. As noted in Lemma~\ref{l:connexite}, the boundary  of   $P_n(x)$  intersects $U$ along two rays 
$R_U(t_n),R_U(t'_n)$ where $t_n, t'_n$ converge to the same value $t$ (since $t_n-t'_n\to 0$).
Denote by  $R_\infty(\zeta_n),R_\infty(\zeta'_n)$ the external rays lying in $\Gamma_n(\theta)$ in ``front of'' $R_U(t_n),R_U(t'_n)$ respectively, {\it i.e.} landing at the same point of $\partial U$.
The sequences $(\zeta_n),(\zeta'_n)$ are monotone and bounded because the puzzle pieces are nested, thus they converge to some limit $\zeta,\zeta'$ respectively.  The rays $R_U(t),R_\infty(\zeta)$ and $R_\infty(\zeta')$ are fixed by $f$ (so they land).   Indeed,  $f(P_{n}(x))=P_{n-1}(x)$ implies that
$f(R_U(t_n))=R_U(t_{n-1})$ and $f(R_U(t'_n))=R_U(t'_{n-1})$ so that $f(R_\infty(\zeta_n))=R_\infty(\zeta_{n-1})$ and $f(R_\infty(\zeta'_n))=R_\infty(\zeta'_{n-1})$.
 Therefore the limit angles $t$ and $\zeta,\zeta'$ are fixed by multiplication by $d$ an $D$ respectively. 
 Let $y,z$ and $z'$ be the landing point of   the rays $R_U(t),R_\infty(\zeta)$ and $R_\infty(\zeta')$ respectively. They ly  in $End(x)$, so if $End(x)=\{x\}$ they coincide.
 
 Assume that  $End(x)\neq\{x\}$. In order to prove that   $x=y=z=z'$, we  look at the "external class" of $f$ on $End(x)$ following an argument of  J.~Kiwi (see~\cite{Kiwi}).  
\vskip 1 em
{\it  $*$ External class $\ol g$ of $f$ on $End(x)$} :

\noindent 
 The set $End(x)$ is a non trivial compact full connected set
 (as the  intersection of a decreasing sequence of  such sets),
 so we can consider a Riemann map $\phi:\overline\C\setminus End(x)\to
\overline \C\setminus \overline \D$.  
First of all, if  $End(x)\cap\partial P_0(x)\neq
\emptyset$
(this happens only in the parabolic case),
we consider  a small enlargement $U_0$ of $P_0(x)$ at
those points such that $U_0\setminus End(x)$ does not contain critical values of $f$.
Let $V_0=\phi(U_0\setminus End(x))$, $U_1=f^{-1}(U_0\setminus End(x))$
and $V_1=\phi(U_1)$. Then, the map $g=\phi\circ f \circ
\phi^{-1}$ is well defined from $V_1$ to $V_0$,  since there is no
pre-image of $End(x)$ in $P_0(x)$ other than $End(x)$.
 Applying Schwarz reflection principle on $V_1$ and
$V_0$, we get neighborhoods   $\tilde V_1$ and $\tilde V_0$ of $\Sun$
and a map $\tilde g :\tilde V_1\to \tilde V_0$ such that  $\tilde
g_{|_{V_1}}=g$. Since $\tilde g$ is a holomorphic covering that preserves
$\Sun$ and each side, 
it has  no critical point on $\Sun$. Therefore, the map  $\ol g=\tilde g_{|_{\Sun}}$ is a covering of $\Sun$, it is called the external class of $f$ on $End(x)$.
\vskip 1 em
{\it  $*$ Fixed points of $\ol g $} :

\begin{lem}
The fixed points of $\ol g$ are weakly repelling, {\it i.e.}  if $p$ is a fixed point of $\ol g$ on $\Sun$ then for any  $z\neq p$ which is close enough to $p$, $| \ol g(z)-p|_{\Sun}>|z-p|_{\Sun}$. \end{lem}
\proof  (\cite{Kiwi,Ro}) Assume, in order to get a contradiction, 
that $p$ is attracting  on one side (at least). 
Then $p$ is  an attracting or a parabolic fixed point for  the holomorphic map $\tilde g$. Therefore  there exists an arc 
$\alpha\subset V_1$ from $q\neq r$ two  points of  $ \S^1$  (one is equal to $p$ in the parabolic case) which bounds an open set 
$\Omega_1$ in $\C\setminus \overline \D$ such that $\tilde g(\Omega_1)\subset \Omega_1$ (it is in the attracting domain of $p$). 
Let  $\Omega=\phi^{-1}(\Omega_1)$. Therefore, $f(\Omega)\subset \Omega$, so that the family $(f^n)$ is normal on $\Omega$. This gives  a ``half-neighbourhood'' of points of $J$ on which $f$ is normal. We prove now that this is impossible for polynomials. Let    $\wt\Omega$  be the Fatou component containing  $\Omega$\,; it is  bounded and fixed by $f$, moreover, by the  Denjoy-Wolf Theorem, $\ol{\wt\Omega}$ contains a fixed point which  attracts every point of $\wt\Omega$.
Moreover $\partial \Omega \cap End(x)\subset\partial \wt \Omega$ contains  more than one point, so there is a cross cut $c$ of $End(x)$ in $\tilde \Omega$. Indeed, the map $\phi^{-1}$ admits  limit points at 
almost every $\theta \in(p,q)$  by Fatou's Theorem and these limit points are not all equal 
(see for  instance Koebe Lemma page 31 of \cite{Go}). So we can use two such part of rays to construct an arc in $\tilde \Omega$ whose end points are on $End(x)$. 
This  cross cut  $c$  bounds a domain $W$ such that  every point  of $\partial W\cap\partial \wt \Omega$ is in $End(x)$. 
Note that $W$ is the trace in $\wt \Omega$ of some open set $W'$ intersecting $J$ ({\it i.e.} $W=W'\cap\wt\Omega$). There is some $i>0$ such that $f^i(W')\supset J$ and therefore $f^i(\partial W\cap\partial \wt\Omega)=\partial \wt \Omega$.  Hence, since  $\partial W\cap\partial \wt\Omega\subset End(x)$ and $f^i( End(x))\subset End(x)$ it follows that  $\partial \wt \Omega\subset End(x)$. This is not possible since $End(x)$  is a full compact connected set  and $ \wt\Omega$ is a topological disk.
One can also notice that  the map $f$ on $\wt \Omega$ is conjugate to a Blaschke product $b:\D\to\D$ with Julia set $\partial \D$ (so that at most one point of $\tilde \Omega$ is not in $\cup_if^i(W)$) .
\cqfd

\noindent 
\begin{claim}Let $u$ be a fixed point of $f$ in $End(x)$ which is accessible from $\ol \C\setminus End(x)$ by an access $\delta$ which is fixed by $f$. Then $\phi(\delta)$  is fixed by $g$ and lands at a point $v:=\phi((u,\delta))$ of $\S^1$ which is also fixed by $g$. \end{claim}
\proof  This follows from the classical theory of Riemann maps (see~\cite{CG,Po}).
\cqfd

\vskip 1 em
{\it  $*$ Some fixed points   of $\ol g $} :
\vskip 0.5 em

\noindent  The points $y,z,z'$ are fixed  points of  $f$  and have respective  accesses   $R_U(t),R_\infty(\zeta),R_\infty(\zeta')$  which are   fixed by $f$.
We can choose $\phi$  up to composing with some rotation such that the landing point $\phi((y,R_U(t)))=1$. Denote by  $u,u'$ the images of  the landing point of $\phi((z,R_\infty(\zeta ))),\phi((z',R_\infty(\zeta')))$ respectively.
 
 Note that the point $y$ is also accessible by an external  fixed ray, say $R_\infty(\eta)$ (since it is repelling or parabolic), and 
 the point $\phi((y,R_\infty(\eta)))=1$.

\vskip 1 em
{\it  $*$ Localization of the fixed points of $\ol g $} :

\noindent 
\begin{claim}Between two fixed points of $\ol g$ there is  a  strict preimage of $1$ by $\ol g$.
\end{claim}
\proof Consider a lift $\hat g$ of $\ol g$ from $\R$ to $\R$ such that $\hat g(0)=0$. It is a strictly monotone map that sends $[0,1]$ to $[0,r]$ for some $r$. Since the fixed points of $\hat g$ are weakly repelling, the  graph of $\hat g$ crosses the line $D_k$ of equation $y=x+k$ for $k\in\{0,\cdots, r-1\}$ from below to above. Hence, between the crossing with $D_k$ and $D_{k+1}$,  the graph of $\hat g$ crosses the line of equation $y=k+1$. The result follows.
\cqfd

\begin{corollary}The rays $R_\infty(\eta), R_\infty(\zeta), R_\infty(\zeta')$ land at the same point.
\end{corollary}
 \proof Assume in order to get a contradiction  that $u\neq 1$
 (recall that $u:=\phi((z,R_\infty(\zeta )))$ and that  $ 1:=\phi((y,R_U(t)))$).  Since the curve $\phi(R_\infty(\eta))$ is an  access to $1$ from $\C\setminus \ol \D$, there is a preimage by $\tilde g$ of $\phi(R_\infty(\eta))$ 
  in each connected component of $\C\setminus (\ol \D\cup \phi(R_\infty(\eta)) \cup \phi(R_\infty(\zeta)))$, by the claim.
  Therefore,   for the open set $U_1\setminus
End(x)$  the following holds :  there are preimages $R_\infty(\eta'),R_\infty(\eta'')$ of $R_\infty(\eta)$  by $f$ 
  inside each connected component of $\C\setminus (End(x)\cup R_\infty(\eta) \cup R_\infty(\zeta))$.
Assume  that $\zeta<\zeta'$ so that  $\zeta_n<\zeta<\zeta'<\zeta'_n$ because the puzzle pieces are nested (the other case is identical).
Since these rays enter every
puzzle piece $P_n(x)$, at least  one of the angle
$\eta',\eta''$ belongs to the intervals $(\zeta_n,\zeta)$ for every $n\ge 0$. 
Therefore, it is equal to $\zeta$ 
but this is impossible since it is strictly pre-fixed. Thus $u=1$. By the same reason,  $u'=1$, so that $R_\infty(\eta)$, $R_\infty(\zeta)$ and
$R_\infty(\zeta')$ land at the same point $z$.

\vskip 1 em
{\it  $*$ The curve $  R_\infty(\zeta)\cup  R_\infty(\zeta')\cup\{z\}$ separates $\ol U$ from $End(x)\setminus\{x\}$} :
\vskip 0.5 em

\noindent 
Let $W_0$ (resp. $W_1$) be the union of the connected components of  $P_0(x)\setminus(R_\infty(\zeta)\cup R_\infty(\zeta')\cup\{z\})$ (resp.  of  $P_1(x)\setminus(R_\infty(\zeta)\cup R_\infty(\zeta')\cup\{z\})$) intersecting $U$.
 Assume to get a contradiction that
$I:=End(x)\cap W_1$ is not empty. Either the map $f:W_1\to W_0$ is an homeomorphism or $I$ contains a critical point (there is no critical point in $P_0(x)\setminus End(x)$).
In the second case, there is  a preimage $R_\infty(\eta''')$ of $R_\infty(\eta)$
landing at a preimage of $z$ in $I$ (since $f(I)=I$) and the angle $\eta'''$ belongs
  to one of the intervals $(\zeta_n,\zeta)\cup (\zeta',\zeta'_n)$ whose diameters tend
to $0$, thus $\zeta=\eta'''$ or  $\zeta'=\eta'''$, which gives the contradiction. 
In  the first case, $f_0=f_{\vert_{W_1}}^{-1}:W_0\to W_1$ is a conformal map and
since $W_1\subset W_0$, by Denjoy-Wolff's Theorem (see~\cite{Steinmetz}),  there exists a fixed point $x$ in $\ol W_0$ to which $f_0^n$ converges uniformly on every compact set of $W_0$. If $x$ is a parabolic point, let $P$ be an open  repelling petal  near $x$. By definition for any  some small $\epsilon>0$,  there  exists some $N>0$, such that  $f_0^n(P)$ is in the $\epsilon$ neighbourhood of $x$ for any $n>N$. 
Now, since $f_0^n$ converges uniformly  to $x$  on the compact set $I\setminus P$, there exists some $M>0$ so that  for $n>M$, $f_0^M(I\setminus P)$ is in the $\epsilon $ neighbourhood of $x$. This contradicts the fact that $f_0(I)=I$ since $f_0^n(I)=f_0^n(P)\cup f_0^n(I\setminus P)$ is in the $\epsilon$ neighbourhood of $x$.
(If $x$ is repelling the argument is even easier.)
Hence the curve  $\ol R_\infty(\zeta)\cup \ol R_\infty(\zeta')$ separates $\overline U$ from $End(x)\setminus\{x\}$.

This achieves the proof of Proposition~\ref{p:kiwi}.\cqfd

\begin{corollary}\label{c:pointfixe}Let $y\in \partial U$ be a point on a preimage $\Gamma_N$ of the graph $\Gamma$ defining the puzzle. Let $(P_n)$ 
be a nest of puzzle pieces with $y$ as a    common boundary point.
Then $(\cap \overline {P_n})\cap \partial U =\{y\}$.
\end{corollary}
\proof  For every $n\ge 0$, we have  $ \overline {P_n}\cap \partial U\subset I_n(y)$, where $I_n(y)$ is  given by Definition~\ref{d:changementdegraphe} for the graph $\Gamma$. Now, since  $y$ is on $\Gamma_n\cap \partial U$, it  is an (eventually)  periodic point and  Corollary~\ref{c:kiwi} implies that $\cap I_n(y)=\{y\}$ (and this holds for any graph   by Lemma~\ref{l:changementdegraphe}). Therefore  $(\cap \ol{P_n})\cap\partial U=\{y\}$.\cqfd

\section{The classical techniques}\label{s:techniques}
 Now in the rest of the paper we consider points  $x\in \partial U$  that are {\it not  eventually periodic}. We will consider graphs $\Gamma$ such that $x$ is   not  on a preimage of  $\Gamma$, so that  the puzzle pieces $ P_n(x)$ are well defined.
  We will prove  that $\bigcap_{n\in \N} \ol P_n(x)$reduces to  $\{x\}$.
  \subsection{Modulus techniques}
One  classical way to obtain that $End(x)=\{x\}$ is to consider the modulus of  the annuli $A_{n_k}:= P_{n_k}
\setminus \ol{P_{n_{k+1}}}$ where $P_{n_k}$ is some  subsequence of the puzzle pieces containing $x$ (not necessarily with consecutive puzzle pieces) satisfying   $\ol{P_{n_{k+1}}} \subset{P_{n_{k}}}$. 

\begin{definition}Any open annulus $A$ is conformally equivalent to $D(0,R)\setminus \ol{D(0,1)}$ for some $R>1$. Then  its modulus is $\displaystyle \mod(A)=\frac{\log R}{2\pi}$. \end{definition}
The following property is classical (see~\cite{Mlc} or Lemma~1.17 of~\cite{Ro1})
\begin{lem} Let   $\Gamma$ be a graph defining a puzzle and let $x$ be any point which is not on the iterated pre-images of $\Gamma$\,; for  some  $n> k$ integers, denote by $D$ the degree of $f^{n-k}:P_n(x) \to  P_k(y)$, where $y=f^{n-k}(x)$. If   there exists  some $r>0$ such that  $\ol{P_{k+r}(y)}\subset  P_k(y)$ then $\ol{P_{n+r}(x)}\subset P_n(x)$ and $\displaystyle \mod (P_{n}(x)\setminus \ol{P_{n+r}(x)})\ge 1/D \mod (P_{k}(y)\setminus \ol{P_{k+r}(y)})$. The equality holds if there is no critical points of  $f^{n-k}$ in $P_{n}(x)\setminus \ol{P_{n+r}(x)}$.
\end{lem}
The following Lemma will be used several times in the text. 
\begin{lem}\label{l:liminfcrit}If $\liminf \mod(A_{n_k})>0$ then $End(x)=\{x\}$.
\end{lem}
\proof The annuli $A_{n_k}$ are disjoint and essential in $P_0(x)\setminus End(x)$.  Up to extracting a sub-sequence, there exist $\epsilon>0$ such that for any $k\ge 0$  $\mod(A_{n_k})>\epsilon$. Therefore, Gr\"otzsch  inequality implies that $\mod( P_0(x)\setminus End(x))\ge \sum_{k\in \N} \mod(A_{n_k})=\infty$.
Since $P_0(x)$ is bounded in $\C$, it follows that $End(x)=\{x\}$. See~\cite{Ah}. 
 \cqfd
In order to apply this Lemma, we need  first  to construct  annuli between puzzles pieces, {\it i.e.} there is no intersection between their boundaries. Then   to control their moduli
we will use the dynamics and look after long iterates of $f$ with controlled degree.\subsection{Finding non degenerate annuli}
\begin{lem}\label{nondegfix} Let  $y\in \partial U$ be a fixed point. Any graph  $\Gamma$, that does not contain fixed rays, defines a puzzle such that\,:
\begin{itemize}\item $y$ is not on the graph and  $\ol {P_1}(y)\subset P_0(y)$ or\,;
\item $y$   is parabolic and is  on the graph $\Gamma$\,; denote by  $ {P_1}(y),P_0(y)$   the puzzle pieces containing $y$ in their closure of depth $1$ and $0$ respectively, then    their closure meets only  at $y$.
\end{itemize}
\end{lem}
\proof If $y$ is a repelling fixed point, it is the landing point of a fixed ray of $U$. We can suppose (up to changing the coordinate in $U$) that it is the ray of angle $0$. Then by 
Lemma~\ref{l:aboutdistray}  the point $y$ is not on the graph $\Gamma$. Now order the cycle of internal rays in $U$ defining the graph  by $\theta_1<\cdots<\theta_r$, assuming that $0<\theta_1<\theta_r<1$.  Then it is immediate that the puzzle  piece $P_1(y)$ (resp. $P_0(y)$) is bounded in $U$ by the rays of angle $\theta_1/2$ (resp. $\theta_1)$ and $\theta_r/2+1/2$ (resp. $\theta_r$). Since $\theta_1/2<\theta_1$ and $\theta_r/2+1/2\ >\theta_r$, the boundaries of $P_0(y)$ and $P_1(y)$ do not touch in $\ol U$. Note that the boundary of $P_0(y)$ consists only in these two internal rays together with two external rays landing at the corresponding points of $\partial U$ and internal/external equipotentials. Then it is clear that $\ol {P_1(y)}\subset P_0(y)$.

If $y$ is a parabolic fixed point but its immediate basin is not $U$, then $y$ is also the landing point of a fixed ray in $U$ and the proof goes exactly as before.

If $y$ is a parabolic fixed point and $U$ is its immediate basin. Then by definition of the graphs, $y$ belongs to any  graph $\Gamma$.    Since there is no external ray landing at $y$ (it would be fixed since $U$ is fixed), there is only one puzzle piece containing $y$ in its closure. And, as before the puzzle piece $P_0(y)$ is bounded by a parabolic equipotential, two parabolic rays in $U$ and two external rays with an external equipotential. Therefore, the closure of 
$P_1(y)$ and of $P_0(y)$ can meet only at $y$.
\cqfd

\begin{lem}\label{l:nondeg}
Let  $x\in\partial U$ and let  $P$  be any puzzle piece  containing    infinitely many points  in  the orbit of $x$ that we denote by  $\mathcal O_{(k)}^+(x)=\{f^{k_n}(x)\mid n\ge 0\}$.
 If there is no periodic point in the accumulation set  of $\mathcal  O_{(k)}^+(x)$ then,  
 there exists  a puzzle piece  $Q$  satisfying  $\ol Q\subset P$ which contains infinitely many points of $\mathcal O_{(k)}^+(x)$.
 \end{lem}
\proof We assume, in order  to get a contradiction,   that for any $r>r_0$ and  for any $z$ such that  $P_r(z)\subset P$ and $\mathcal O_{(k)}^+(x)\cap P_{r}(z)$ is infinite,  the boundaries of the  pieces intersect\,:  $\partial P\cap \partial P_{r}(z)\neq\emptyset$.  Since the graph is of the form of $\Gamma(\theta)$ (definition~\ref{d:graph}), the puzzle pieces  can only intersect  either at the parabolic point or at the landing point of finitely many rays   on $\partial P$, we deduce that  for infinitely many values of $r$
 they will intersect at the same point $v$. This point has to be  an eventually periodic point  of $\partial U$. Thus we get a nested sequence $(P_n)$ of puzzle pieces containing a common point  $v$ in their boundary   and by Corollary~\ref{c:pointfixe} we know that 
 $\cap (\ol{P_n}\cap\partial U)=\{v\}$.
Since, for any $r\ge 0$, $ \ol{P_{r}}(z)\cap \partial U$ contains infinitely many points of   $\mathcal  O_{(k)}^+(x)$ and  $\cap_{r\ge r_0} \ol{P_{r}}(z)\cap \partial U=\{v\}$,  $v$ is in the accumulation set of $\mathcal  O_{(k)}^+(x)$. This contradiction with  the assumption implies  the Lemma.\cqfd

\subsection{Distortion property}
 The property for an annulus to be non degenerate  is difficult to obtain in the parabolic case where puzzle pieces touch at the pre-images of the parabolic point.
Therefore,  we have to enlarge the pieces to obtain a non degenerate annulus.
Doing this we might lose  the property that the annuli are disjoint. 
So we use the following  distortion Lemma (see~\cite{Yin,YinZhai}). Note that we still need to have some non degenerate annulus and also bounded degrees.
\begin{lem}\label{l:distortionlemma}Let $U,V, U_n,V_n$  be topological disks  with $\ol V\subset U$ and  such that for some $k_n$,  $f^{k_n}(U_n)= U$  and $f^{k_n}(V_n)=V$. 
Moreover, if  $V_n$ contains some point $z$ of $J(f)$ and if  the degree of the maps $f^{k_n}\from U_n \to U$ is  bounded independently of $n$, then the diameter of $V_n$ tends to $0$.\end{lem}
\proof  Assume by contradiction that the diameter of $V_n$ does not  tend to $0$. Up to extracting a subsequence, there exists $L>0$ such that $diam(V_n)\ge L$.
Denote by $Shape(U,z)$ the ratio $\displaystyle \frac{\max_{t\in\partial U} d(t,z)}{\min_{t\in\partial U} d(t,z)}=\frac{\max_{t\in\partial U} d(t,z)}{d(z,\partial U)}$. We  will prove below that the shape $Shape(U_n,z)$ is bounded  above by some constant $S>0$ independent of $n$. In particular, the internal diameter is bounded below by a positive constant independent of $n$\,:
Since  the maximum $\max_{t\in \ol{U_n}} d(z,t)$ is reached by some  $t$ in $\partial U_n$, one has  that $diam(U_n)\le 2 d(z,\partial U_n)$. Therefore, $\min_{t\in\partial U_n} d(z,t)=\max_{t\in\partial U_n} d(t,z)/Shape(U_n,z)\ge L/(2S)$. Hence there is a small disk $D$ around $z$  contained in every $U_n$. This contradicts 
with $J(f)\subset f^k(D)$ for large values of $k$ since $z$ is in $J(f)$.

\vskip 0.2em
To prove that the shape is bounded above, we use the following property on the shape (see  Lemma 2 of \cite{YinZhai})\;:

{\bf Lemma\,:} {\it Let $\Delta$ be the unit disk, $ U, \tilde U, V$ and $ \tilde V $  be topological disks and $g:(\Delta, U, \tilde U) \to (\Delta, V, \tilde V )$ be a holomorphic proper map of
degree $d$ with  $0 \in \tilde U \subset U \subset \Delta$, $ 0 \in \tilde V \subset V \subset\Delta$. Suppose that 
$ deg(g |_{\tilde U}) = deg(g |_{U}) = deg(g |_{\Delta}) = d \ge 2$. Then there exists a constant $K = K(m, d) > 0$, where 
$ mod(V \setminus \tilde V ) \ge  m > 0$,
such that
$Shape(U, 0) \le K á Shape(V, 0)^\frac1d$ .}

The proof of this Lemma uses the classical Koebe distortion Theorem for univalent function and the  Gr\"otzsch Theorem.

\subsection{Bounded degree, property $(\star)$ and   the successors.}
To get control on the moduli or to apply the distortion Lemma, we need to have maps of bounded degree.
\begin{definition}Let  $(\star)$ be the following property of $x$ relative to a graph $\Gamma$\,:    
\begin{center}$(\star)$\,: There is a sequence $(k_n)_{n\ge 1}$ tending to $\infty$, a point $z$ and a level $k_0$ such that 
 the  degree of the maps $f^{k_n}:P_{k_n+k_0}(x)\to P_{k_0}(z)$ is bounded  independently of $n$.
 \end{center}
\end{definition}
One simple way to get a bounded  degree  is to consider {\bf the first entrance} or the  {\bf the first return}  to a puzzle piece\;: 
\begin{definition}The  {\it the first entrance time} (resp.  {\it the first return time})  of  a point $z$  in puzzle piece $P$, is the minimal  $r\ge 0$  (resp. $r\ge 1$) such that 
 $f^r(z)\in P$. We call  {\it the first entrance in $P$} (resp.  {\it the first return in $P$})  the point  $f^r(z)$.
 \end{definition}
\begin{notation}\label{n:numbercritend}Let $b$ denote the number of distinct critical  ends, {\it i.e.}  ends containing at least one critical point (see Notation~\ref{n:critend}). 
Let $\delta$ be the maximum of the degree of $f$ over all the critical  ends (see Remark~\ref{r:enddegree}).
\end{notation}

\begin{lem}\label{l:firstentrance}Let $r\ge 0$ be the first entrance time (resp. the first return time) in some puzzle piece $P_n$ (of depth $n$) of a point $z$, then the degree of $f^r:P_{n+r}(z)\to P_n$ is bounded by $\delta^b$ (resp.  $\delta^{b+1}$).
\end{lem}

\proof  The Lemma  follows from the fact that the sequence of puzzle pieces $\{f^i(P_{n+r}(z)) \mid 0\le  i\le r\}$ meets each of the $b$ critical  ends at most once. Indeed, assume in order  to get a  contradiction, that a critical point $c$ belongs to $f^i(P_{n+r}(z))=P_{n+r-i}(f^i(z))$ and to $f^j(P_{n+r}(z))=P_{n+r-j}(f^j(z))$ for $0\le i<j\le r$. Then the point $f^i(z)$ which is in $P_{n+r-i}(c)\subset P_{n+r-j}(c)$ is mapped by $f^{r-j}$ in the piece $P_n=f^{r-j}(P_{n+r-j}(c))$. This  contradicts the fact that $r$ is the first entrance time of $z$ in $P_n$  since $0\le r-j+i<r$. The Lemma  follows in this case. If $r$ is the first return time we get a contradiction if $ r-j+i\ge 1$ which happen exactly when $j=r$ and $i=0$. In this case, on see only  one critical end twice so the degree is bounded by  $\delta^{b+1}$.
\cqfd

One particular case where the bounded degree property  appears is when we look at the map from a {\bf successors} of a piece $P$ to the piece $P$. They are, in the case of several critical points, the generalisation of  the notion of {\it children}    introduced by Branner and Hubbard. Recall that\,: \begin{itemize}
\item
 a {\it piece $P_{n+k}(c)$ is called a child of $P_n(c)$}  if   the map $f^{k-1}:P_{n+k-1}(f(c))\to  P_n(c)$ is a homeomorphism where $c$ is a critical point\,;
 \item in the case of  several critical points,     it can be generalized by\,: 

{\it a piece $P_{n+k}(c')$ is called a child of $P_n(c)$} if  $f^{k-1}:P_{n+k-1}(f(c'))\to  P_n(c)$ is a homeomorphism where $c,c'$ are critical points. 
\end{itemize}
 To define a nest around a given critical point, we need to come back  to the same  critical point.
\begin{definition}\label{d:successor}Let $c$ be a critical point. The piece  $P_{n+k}(c)$ is said to be a {\it successor} of $P_n(c)$ if $f^k(P_{n+k}(c))=P_{n}(c)$ and each critical point meets at most twice  the set of pieces $\{f^i(P_{n+k}(c)),\ 0\le i \le k\}$.
\end{definition}
Note that   $P_{n+k}(c)$ and  $P_{n}(c)$ contain $c$, therefore the set of pieces  $\{f^i(P_{n+k}(c)),\ 0\le i \le k\}$ contains $c$  at least  two times. For this reason, we cannot impose less than twice in the definition. 
 \begin{remark}\label{r:successor}If $P_{n+k}(c)$ is a successor of $P_n(c)$, then the map $f^k:P_{n+k}(c) \to P_{n}(c)$ has a degree  less than $\delta^{2b}$.
\end{remark}
\begin{corollary}\label{c:nonpersist} Let  $c$  be a  critical point. If there exists $n_0\ge 0$ such that the piece $P_{n_0}(c)$
 has infinitely many successors, then $c$  has property  $(\star)$.\end{corollary}
\proof If  $P_{n_0+k_n}(c)$ is an infinite  sequence of successors, it follows from  the last remark.
\cqfd

\section{Combinatorial accumulations.}

To find long iterates starting from  any puzzle piece   $P_n(x)$ with bounded degree, we need to consider  the critical points which appear in the orbit  of  $P_n(x)$. More generally we will have to consider   how the critical points  accumulate themselves.
\begin{definition}\label{d:accumulation} For a given graph $\Gamma$, defining a puzzle, we say that $z'$ is in the  {\it  combinatorial accumulation  } of $z$ and  we note  $z'\in \wc(z)$ if\,:
$$\forall n\ge 0\quad \exists k> 0\quad \hbox{ such that } f^k(z)\in P_n(z').$$
\end{definition} 
\begin{remark}\label{r:accumulation} The relation is clearly transitive. Moreover\,:
 \begin{itemize}
\item The notion of combinatorial accumulation depends on the graph $\Gamma$, at least by the fact that a point $z'\in \Gamma$   cannot be in $\wc(z)$\,;
\item  Any iterate  $ f^k(z)$ 
of $z$  (with $ k>0 $) is in $\omega_{comb}(z)$\,;
\item
This  notion coincides with the standard notion of $\omega$-limit set $\omega(z)$ but for the topology  generated by the set of all  the puzzle pieces, except as we just noticed, that it contains the orbit of the point $z$ and that it does not contains the points of  $\omega(z)$ which are on the iterated preimages of $\Gamma$\,;
\item In particular, $z'\in \omega(z)$  and $f^n(z')\notin \Gamma$ for all $n\ge 0$ implies that $z'\in \wc( z)$\,;
\item
For the converse, if $z'\in \wc( z)$ and  if $End(z')=\{z'\}$ (resp.   $End(z')\cap\partial U=\{z'\}$), then either $z'$ is  in the orbit of $z$ or $z'\in \omega(z)$.
\end{itemize}
\end{remark} 
\begin{remark}\label{r:accumulatingperiodic}
Let $x\in \partial U$, the combinatorial accumulation does not depend on the graph (see Lemma~\ref{l:changementdegraphe}) in the following sense\,:
$$y\in \wc(x)\iff \forall n\ge 0 \ \exists k\ge 0\hbox{ such that }f^k(x)\in I_n(y).$$ 
Moreover,   if   $y$ is eventually periodic by Corollary~\ref{c:kiwi} implies that   $$y\in \wc(x)\iff y\in \omega(x) \hbox { or }
y=f^k (x) \hbox{ for some } k\ge 0.$$ 
\end{remark} 
From the definition (Definition~\ref{d:accumulation}) it follows that one can always make the following assumption,  up to replacing the graph $\Gamma$ defining the puzzle  by some of its iterated  pre-image $\Gamma_N$ (note that the combinatorial accumulation set is the same for $\Gamma$ and $\Gamma_N$)\,:
\begin{assumption}\label{a:nocrit}Let $Y$ be a finite set of points containing the critical points.  We can assume that
 for any $y,y'\in Y$,  if $y\notin \wc(y')$ then $y\notin P_0(f^k(y'))$ for all $ k\ge 0$ {\rm(}equivalently $y\in \wc(y')$ or $f^k(y')\notin P_0(y)$ for all $ k\ge 0${\rm)}. \end{assumption}
 \proof By definition,   for any  two points $y,y'$,  either $y\in \wc (y')
$ or there exists $N\ge 0$ such that 
$y\notin P_n(f^k(y'))$ for any $n\ge N$ and any $k\ge 0$. If $y$ is on some iterated preimage of $\Gamma$ we  fall in the second case. 
Since  the set $Y$  is  finite,  taking the  largest $N$ when  $y, y'$ belong to $Y$
and taking the graph defined by  the  $N$-th preimage of $\Gamma$, the puzzle pieces  avoided become of level~$0$.\cqfd
%
\subsection{Periodic point in the accumulation}
In this subsection we take a point $x\in \partial U $ which is not eventually periodic and which contains a periodic point $y$ in its $\omega$-limit set. We can assume that $y$ is fixed up to taking an iterate of $f$, since the Julia set (and also the boundary of the basin $U$) will be identical for $f$ and for the iterate.  
\begin{proposition}\label{p:accufix} Let $x\in \partial U$. Assume that $x$ is neither  periodic  nor  eventually periodic and that  it contains  a fixed point $y$ in its $\omega$-limit set. 
Then, there exists a graph $\Gamma$ defining a puzzle such that $End(x)\cap\partial U=\{x\}$. In particular, $\partial U$ is locally connected at $x$. 
\end{proposition}
\begin{corollary}\label{c:accufix} Let $x\in \partial U$. Assume that $x$ is neither  periodic  nor  eventually periodic and that  it contains  a periodic  point  $z$ in its $\omega$-limit set. 
Then  $\partial U$ is locally connected at $x$.
\end{corollary}
\proof The corollary follows directly from the proposition by taking some iterate $g$  of $f$ and for   $y$ the adapted iterate forward image of $y$ by  $f$. \cqfd

\vskip 1em
\noindent {\it Proof of Proposition~\ref{p:accufix}.}  
Since $y\in \partial U$,  for  any graph $\Gamma$ that do not contain fixed rays,  we can define the sequence of puzzle pieces $(P_n(y))$ (that contains $y$ in its closure).  Moreover, for this graph $\ol{P_1(y)}\subset P_0(y)$ if $y$ is repelling and   $\ol{P_1(y)}\cap \ol{ P_0(y)}=\{y\}$ if $y$ is parabolic on $\Gamma$
 (by Lemma ~\ref{nondegfix}). 
 Note that if we have to take an iterated preimage of $\Gamma$ in order to satisfy Assumption~1 or~2,  we will keep this property\,:  the preimage  (by $f^r$) containing $y$ of the puzzle pieces $P_1(y)$ and $ P_0(y)$ are puzzle pieces $P_{r+1}(y)$ and $ P_r(y)$  that satisfies $\ol{P_{r+1}(y)}\subset P_r(y)$ in the first case and $\ol{P_{r+1}(y)}\cap \ol{ P_r(y)}$ reduces to some preimages of $y$ by  $f^r$ in the second case.

 For any $n\ge 0$, let $r_n$ be the first entrance time of $x$  into $P_n(y)$. 
 The sequence $r_n$ is not eventually constant. Otherwise, we would have   $f^r(x)\in End(y)$ for  $r=r_n$ with $n\ge N_0$, but $f^r(x)\in \partial U$ and 
 $End(y)\cap \partial U=\{y\}$ (Corollary~\ref{c:kiwi}) would imply that $f^r(x)=y$ (contradiction).

Then,  no iterate of $x$ is in $End(y)$, so  for  $n\ge 0$, there exists $s> n$ such that $f^{r_n}(x)\notin P_{s}(y)$.  Choose $s_n$ the  minimal  $s>n$ with this property. Now,  for any $i\ge 0$, $P_{i+1}(y)$ is the only preimage of $P_i(y)$ in $P_0(y)$  since  there is no critical point in $P_0(y)\setminus End(y)$ (by Assumption~\ref{a:crit}). Therefore,   $f^j( P_{s_n}(f^{r_n}(x)))\subset  P_{{s_n}-j-1}(y)\setminus P_{{s_n}-j}(y)$ for $j\le {s_n}-1$ and $f^{{s_n}}:P_{{s_n}}(f^{r_n}(x))\to P_0(f^{r_n+{s_n}}(x))$ is a  homeomorphism\,; but
$f^{r_n+s_n}(x)\notin  P_0(y)$ (by the same reason as before).

We should return to $P_1(y)$ but with bounded degree, and  also bounded degree on the puzzle that is mapped to $P_0(y)$. 
\vskip 1em
 1) {\it We consider first the case where $y$ is a parabolic point}. 
Taking the first return time of  $f^{r_n+{s_n}}(x)$ to $P_0(y)$ 
give now a sequence  $m_n\ge n$ such that  the map 
$f^{m_n}:P_{m_n}(x)\to P_0(y)$ has a  degree bounded by $\delta^{2b}$. Moreover, 
  $y$ is not in the closure of 
 $f^{m_n-1}(P_{m_n}(x))$ otherwise  $f^{m_n-1}(P_{m_n}(x))$ would already be in  $P_0(y)$. Therefore,  $f^{m_n-1}(P_{m_n}(x))$   contains a preimage $y'$ of $y$ which lies in  a preimage $U'$ of $U$. We  enlarge this puzzle piece
$P_1(y')$ to  get $\tilde P_1(y')$ (an open set) which is the  union of $P_1(y')$ with a small neighborhood of $y'$ in $U'$ that avoids the orbits of the critical points. This is possible since  $U'$  is a pre-image of the basin $U$, so the orbits of the critical points intersect under a finite set $U'$ where we enlarge the puzzle piece.
Let $\widetilde{P_{m_n}}(x)$ denote the iterated preimage  of  $\widetilde{P_1}(y')$ by $f^{m_n}$ (it is an enlargment of $P_{m_n}(x)$). The degree of $f^{m_n-1}:\widetilde{P_{m_n}}(x)\to \widetilde{P_1}(y')$ is exactly the same as the degree of $f^{m_n-1}:{P_{m_n}}(x)\to {P_1}(y')$.  Lemma~\ref{l:distortionlemma} apply then,  since we have a non degenerate annulus between $\widetilde{P_1}(y')$ and ${P_2}(y')$. 
The local connectivity follows then from Lemma~\ref{l:connexite}.

\vskip 1em 2){\it  We consider now the case of a repelling fixed point $y$.}  
Here we can not avoid the post-critical set that maybe accumulates on $y$. So we look at the critical points that appear in the orbit   of the puzzle pieces 
 which  are mapped  to $P_0(y)$ (by some  iterate of $f$). More precisely, consider  $A:=
 \{c\in \wc(x)\cap Crit\mid y\in \wc(c)\}$. Let $c$ be a critical point such that  some iterate $f^k(x)\in P_n(c)$ for some $n\ge 0$, then if  for some $0\le i\le n$ a puzzle piece  $P_j(f^i(c))$ contains $y$ (for some $j\ge 0$) then $c\in A$ (by Assumption~\ref{a:crit} and ~\ref{a:nocrit}). 
  
  Assume first that $A=\emptyset$. Let  $t_n$ is  the first  entrance  time of $f^{r_n+{s_n}}(x)$ in $P_1(y)$. There is no critical point  in  $f^{i}(P_{t_n}(f^{r_n+{s_n}}(x)))$ 
  for $0\le i\le t_n$,  otherwise this critical point would belong to $A$.
 Then the map $f^{t_n}: P_{t_n}(f^{r_n+{s_n}}(x))\to P_0(y)$ is a homeomorphism.

  We consider now the case of $A\neq \emptyset$. For  any $c\in A$, we  denote by  $e_c$ the first entrance time    of $c\in A$ in $P_1(y)$ (it exists since $y\in \wc(c)$), and consider   the maximum $N$ of $e_c$ for $c\in A$. For each $c\in A$, we want to consider the first entrance time of $f^{r_n+{s_n}}(x)$ in $P_N(c)$.   Since $c\in \wc(x)$, it is easy to define,  if no image of $x$ is in $End(c)$. Indeed,  for any $k\ge 0$ there exists $i_k$ such that $f^{i_k}(x)\in P_k(c)$, and  if the sequence   $i_k$  does not tend to infinity, one can extract   a sub-sequence that is constant so that  some iterate 
   $f^{i}(x)\in End(c)$. In this case, if $f^i(x)\in End(c)$, we start again with $z =f^i(x)$, if the sequence $j_k>0$, such that  $f^{j_k}(z)\in P_k(c)$  for $k\ge 0$, does not tend to infinite, then it implies that $End(c)$ is periodic (since $End(c)=End(z)$). But since $y\in \wc(c)$ we obtain that for some $j$, $f^j(End(c))=End(y)$ and therefore $f^{i+j}(x)\in End(y)\cap \partial U$ so that  $f^{i+j}(x)=y$ since  $End(y)\cap \partial U=\{y\}$ by Corollary~\ref{c:kiwi}.

   Now let $t_c$ be the first entrance time of  $f^{r_n+{s_n}}(x)$ in $P_N(c)$ and denote by $t_n$ the minimal one for $c\in A$. 
   Let $c_0$ be the critical point such that $f^{r_n+{s_n}}(x)$ meets the first in 
  the level $N$ puzzles pieces, note $e_n=e_{c_0}$ to simplify. Then  the puzzle piece $P_{t_n+e_{n}}(f^{r_n+{s_n}}(x) )$   is mapped by $f^{t_n+e_n}$ to $P_0(y)$ since it 
  is mapped by $f^{t_n}$ to $P_{e_{n}}(c_0)$, which is mapped by $f^{e_n}$ to $P_0(y)$. 
 If ${t_n+e_n}\le N$ the degree of the map   $f^{t_n+e_n}$ is clearly bounded by $\delta^N$. If ${t_n+e_n}> N$, the map $f^{t_n+e_n}$ can decompose into $f^{N+t}$ with $0<t\le t_N$ (since $n_0\le N$). There is no critical point 
  in  $P_{t_n+e_{n}-i}(f^{r_n+{s_n}+i}(x) )$  for $0\le i<t_n$, since any such critical point should be in $A$ and $t_n$ is the first time $f^{r_n+{s_n}}(x) $ meets a critical point in level $N$ puzzle pieces. Therefore the map  $f^{t}: P_{t_n+e_{n}}(f^{r_n+{s_n}}(x) )\to P_N(f^{r_n+{s_n}+t}(x))$ is a homeomorphism and then  the degree of  $f^{N+t}:P_{t_n+e_{n}}(f^{r_n+{s_n}}(x) )\to P_0(y)$ is bounded by $\delta^N$. 
  
  Finally, we get a sequence $u_n=r_n+s_n+t_n+e_n$ tending to infinity such that 
  $f^{u_n}(x)\in P_1(y)$ and the degree of $f^{u_n}:P_{u_n}(x)\in P_0(y)$ is bounded by $D=\delta^{b+N}$. Then Lemma~\ref{l:liminfcrit} applies since $\mod(P_{u_n}(x)\setminus\ol{P_{u_n}(x)})\ge \mod(P_0(y)\setminus\ol{P_1(y)}) /D$.
    \cqfd
\subsection{How to get property ($\star$) and its consequence}
\begin{lem}\label{l:star}
If $y$ has property  $(\star)$ and $ y\in \wc(x)$, then $x$ also has  property  $(\star)$.
\end{lem}
\proof  Suppose $y$ has property  $(\star)$. Then there exist  $z$ and  $(k_n)_{n\in \N^*}$   such that the degree of the maps  $f^{k_n}: P_{k_n+k_0}(y)\to P_{k_0}(z)$ is bounded. 
We look  at the first entrance  of $x$ to the nest of $y$.  For $n> 0$, let   $r_n$ be the first entrance time of $x$ in $ P_{k_n+k_0}(y)$. The  degree of the maps $f^{r_n}: P_{r_n+k_n+k_0}(x)\to P_{k_n+k_0}(y)$ is bounded  by Lemma~\ref{l:firstentrance}. Therefore, the  degree of the maps $f^{r_n+k_n}: P_{r_n+k_n+k_0}(x)\to P_{k_0}(z)$ is bounded 
 and  property $(\star)$ follows for $x$ with the same puzzle piece $P_{k_0}(z) $
 and the sequence $k_n+r_n$.
 \cqfd

\begin{definition}
Let $\omega Crit(z)$ denote the set of critical points  which are  in $\wc(z)$\,:  $\omega Crit(z)=Crit\cap \wc(z)$.
\end{definition}
\begin{lem}\label{l:cassimples}Suppose that $x\in \partial U$ has either one  of the  following properties\,:
\begin{enumerate}
\item   $\omega Crit(x)=\emptyset$\;, 
\item  $\omega Crit(x)\neq\emptyset$  and there exits $c\in \omega Crit(x)$ such that $\omega Crit(c)=\emptyset$\,;
\item    $\omega Crit(x)\neq\emptyset$  and
 there exist $c,c'\in \omega Crit(x)$ such that  $c\notin\omega Crit (c')$ {\rm(}$c$ and $c'$ can be the same{\rm)}.\end{enumerate}
 Then   $x$ satisfies property $(\star)$. 
 
\end{lem}
\proof If $\omega Crit(x)=\emptyset$, the map  $ f^{k-1}: P_{k-1}(f(x))\to P_0(f^k(x))$ is a homeomorphism since no critical point belongs to  $ P_{k-i}(f^{k+i}(x))$ (by Assumption~\ref{a:nocrit}). Therefore the degree of  the map 
$ f^k: P_k(x)\to P_0(f^k(x))$ is bounded by $\delta$ ($x$ can be a critical point). Since there is a finite number of puzzle pieces of level~$0$,  the property $(\star)$ follows.

 \noindent 
For Point 2.,    since  $ \omega Crit(x)\neq \emptyset$, there is a point $c\in  \omega Crit(x)$. If we assume that $\omega Crit(c)=\emptyset$, then  $c$ has property $(\star)$  (by point 1.)  and the statement is a consequence of  Lemma~\ref{l:star}.

\noindent For Point 3., we suppose that Points 1. and 2. are not satisfied\,: $\omega Crit(x)\neq \emptyset$ and     for every  $c\in \omega Crit(x)$, $ \omega  Crit(c)\neq\emptyset$. Note that if there exists $c\in \omega Crit(x)$ such that  $c\notin \omega Crit(c)$, then  for any $c'\in \omega Crit(c)$, we also have  $c\notin \omega Crit(c')$. Therefore in Point 3. we can always take $c\neq c'$.
 Moreover,  for any $r\ge 0$, $f^r(x)\notin End(c')$, otherwise since $c\in \omega Crit(x)$ 
 it would imply that $c\in \omega Crit(c')$. 
 Hence for any $k\ge 0$, if  we  consider the first entrance time $r_k$ of $x$ in $P_n(c')$, there exists a minimal $n_k$ such that $c'\in P_{n_k}(f^{r_k}(x))\setminus  P_{n_k+1}(f^{r_k}(x))$. Now,  by Assumption~\ref{a:nocrit},   for any $r_k\le l\le r_k+n_k$,  $c\notin P_0(f^l(x))$ (since $c\notin \omega Crit(c')$). Then let $l_k$ be the first $l>r_k$ such that $c\in P_0(f^l(x))$ (it exists since $c\in \omega Crit(x)$).  Each critical point meets the set of puzzle pieces 
 $\{f^{i}:P_{l_k-r_k}(f^{r_k}(x)) \mid 0\le i\le l_k-r_k\}$ at most once, otherwise $l_k$ would be smaller. Then the map 
$f^{l_k-r_k}:P_{l_k-r_k}(f^{r_k}(x))\to P_0(c)$ has degree bounded by $\delta^b$.
Now, by Lemma~\ref{l:firstentrance},  the map $f^{r_k}: P_{l_k}(x)\to P_{l_k-r_k}(f^{r_k}(x))$
has bounded degree (since $r_k$ is the first entrance time of $x$ in $P_k(c')$).
Therefore,  the map $f^{l_k}: P_{l_k}(x)\to P_0(c)$ has its degree bounded independently of $k$.
 It follows, by Lemma~\ref{l:star}  that $x$ has property $(\star)$. \cqfd

If $x$  does not satisfy neither 1), 2) nor 3) then the set $\omega Crit(x)$ is non empty and  for any  critical points $c,c'\in \omega Crit(x)$,   $c'\in \omega Crit(c)$ and  $c\in \omega Crit  (c')$ ($c$ and $c'$ can be the same point).

\begin{definition}\label{d:selfrec} A critical point $c\in Crit$ is said  {\it critically self-recurrent}  if $\omega Crit(c)\neq\emptyset$  and every 
 $c'\in\omega Crit (c)$ satisfies  $c\in\omega Crit (c')$. 
 \end{definition}

\begin{remark}\label{r:nonselfrec}
\begin{itemize}\item
Observe that if $c\in Crit$ is critically self recurrent, then $c\in \omega Crit(c)$.\item If there is a point $c\in \omega Crit (x)$ that is not self-recurrent, then $x$ has property $(\star)$. 
\end{itemize}
\end{remark}
\proof The first point follows from  the transitivity property. For the second point, let    $c\in \omega Crit(x)$ be a critical point that is not critically self-recurrent. Then there exists $c'\in \omega Crit(c)$ such that   $c\in \omega Crit(c')$. By the transitivity property, $c'\in \omega Crit(x)$ and Point 3. of Lemma~\ref{l:cassimples} apply.\cqfd

\begin{corollary}\label{c:star} Let $x\in \partial U$  be a non (eventually) periodic that does not accumulate a periodic point. If $x$     satisfies  property $(\star)$  then  $End(x)=\{x\}$. Therefore  $\partial U $ is locally connected at $x$.
\end{corollary}
\proof  Property $(\star)$ gives  a sequence $(k_n)$ tending to $\infty$ and a puzzle piece $P_{k_0}(z)$ such that the degree of the maps  $f^{k_n}:P_{k_n+k_0}(x)\to P_{k_0}(z)$ is bounded independently of  $n>0$.
Assume first that  the accumulation set of $\mathcal O^+_{(k)}(x)$  does not contain  any periodic  point. From Lemma~\ref{l:nondeg}, for some piece $Q$ with   $\ol Q \subset P_{k_0}(z)$, we know that  $Q$
contains infinitely many points of $\{f^{k_n}(x)\mid n>0\}$.
Looking at the component of $Q_n:=(f^{k_n})^{-1}(Q)$  we deduce  that $End(x)=\{x\}$. Indeed, we have  coverings $f^{k_n}$  of degree bounded  by some constant $D$ that sends $Q_n$ to $Q$,  $P_n:=P_{k_n+k_0}(x)$ to $P:= P_{k_0}(z)$ with  $\ol Q_n\subset P_n$, $\ol Q\subset P$, so that $\mod (P_n\setminus Q_n)\ge \frac 1 D  \mod (P\setminus Q)$.
\cqfd

%
\section{The persistently recurrent case}\label{s:persrec}
%
\subsection{Reduction of the set of points}
After the work done in previous sections, we can restrict the set of points we are working with. 
 More precisely, we consider points $x\in \partial U$ that are {\it not eventually periodic }(see section~\ref{s:kiwi}) and  that {\it do not accumulates  on  periodic points} (see Corollary~\ref{c:accufix}).   
 For such a point $x$, if $\omega Crit(x)=\emptyset$ or if there exist $ c,c'\in \omega Crit(x)$ such that $c\notin \omega Crit(c')$   (not self-recurrent) or if there exists 
 $c\notin \omega Crit(x)$  that has infinitely many successors, then $x$ satisfies property $(\star)$ and we get that $End(x)=\{x\}$ by Lemma~\ref{l:cassimples}, 
 Corollary~\ref{c:kiwi},~\ref{c:nonpersist}, \ref{c:accufix} and \ref{c:star}. 
 
Therefore we   consider only points $x\in \partial U$ such that any  critical point $c$ of $\omega Crit(x)$  is 
critically self-recurrent.

\begin{definition}\label{d:perCrit}A    point $c$   is  {\it persistently recurrent } if  it is critically self-recurrent and for any point $c_1\in \omega Crit(c)$, any  $n_0\ge 0$  the puzzle piece $P_{n_0}(c_1)$
 has only finitely many successors.  Denote by  $perCrit(x)$ the set   of points $c$ of $\omega Crit(x)$ which are   persistently recurrent
 \end{definition}

As noticed in Remark~\ref{r:nonselfrec}, if $perCrit(x)\neq \omega Crit(x)$ then $x$ satisfies property $(\star)$. For this reason, {\bf we consider for the rest of the paper only points $x$
 satisfying}\,:\begin{enumerate}
 \item
 $x\in \partial U$ is  not eventually periodic\,; 
 \item $x$ does not accumulate on  periodic points\,;
 \item  and 
 $perCrit(x)=\omega Crit(x)\neq \emptyset$.
 \end{enumerate}

\begin{remark}The critical points that belong to a same end $End(y)$ cannot be separated by puzzle pieces of any depth. Everything goes as if there were only one (at most) critical point in $End(y)$ and this point would have multiplicity equal to $\delta(y)-1$ (see Remark~\ref{r:enddegree}). Our arguments in the proofs will be  as if it were the case. We could have consider the relation on the set of critical points that identifies two critical points if they belong to the same end  and then argue on the equivalences classes but we prefer not to have more notations.   \end{remark} Thus we can make the following assumption without lost of generality.
\begin{assumption}\label{a:3} There is  at most one critical point in each puzzle piece {\rm(}in particular  at level~$0${\rm)}.
\end{assumption}
\subsection{Property of the successors}
In the persistently recurrent case, there are only finitely many  successors, but a 
very  fundamental fact is that there are  always at least two successors. It is the content of the following Lemma whose proof is postponed to section~\ref{s:KSSproof}\,:
\begin{lem}\label{l:2successors}
 Let $c$ be any critical point of $f$.
 \begin{enumerate}
 \item
  If $End(c)$ is not periodic, then each puzzle piece $P_n(c)$ has at least two successors. 
 \item If $c\in \omega Crit(x) $  has an end  $End(c)$ which is  periodic and $x$ is a point of $\partial U$, then $\omega(x)$ contains a periodic point.
 \end{enumerate}
\end{lem}

Note that we do not use in the proof of this Lemma that we have a point $c$ of $perCrit(x)$.
Nevertheless, this Lemma found is utility only in this case (since we do not know if we have enough successors).
\begin{definition}Let $A$ be a puzzle piece, denote by $r(A)$ the smallest  return time of a point of $A$ in~$A$\,:
$$r(A)=\inf_{z\in A}\{k(z)>0\mid f^{k(z)}(z)\in A\}.$$ 
\end{definition}
\begin{remark}\label{r:returntime}\begin{enumerate}
\item If a  point of $A$ returns to $A$, a sub-piece returns by the same iterate  to $A$\,;
\item If $A'\subset A$ are two puzzle pieces, then  $r(A')\ge r(A)$\,;
\item Let $y$ be a point in   a sub-piece $A'$ of $ A$\,;  assume  that for some $k>0$, $f^k(A')=A$  and that $y\notin f^i (A')$ for $0<i< k$, then $r(A')\ge k$.
\end{enumerate}
\end{remark}
\proof If a point comes back to $A'$ it comes back also to $A$ so  $r(A')\ge r(A)$. For the second statement, let $z$ be any point of $A'$.  Note that since $A'$ is a puzzle piece, either  $f^i(A')\cap A'=\emptyset$  or $f^i(A')\supset A'$ (for level reason). Since $y\notin f^i (A')$ for $0<i< k$,    $f^i(A')\cap A'=\emptyset$, so  the points  $f^i (z)$ and $y$  do not belong to the same piece. Hence the first return time of $z$ to $A'$ is necessarily more than $k$.\cqfd
\begin{notation}\label{n:successor}
Denote by $\DD(P)$ the {\it last successor} of  the puzzle piece $P$ and by  $\sigma(P)$   the integer such  that   $f^{\sigma(P)}(\DD(P))=P$.
\end{notation}
\begin{corollary}\label{c:2successors} Let $A$ be a puzzle piece, then $r(\DD(A))\ge\sigma(A)\ge  2 r(A)$.
\end{corollary}
\proof The inequality  $r(\DD(A))\ge\sigma(A)$ follows from previous remark. 

Let $A'$ be the first successor of $A$ and let $k$ be the integer such that $f^k(A')=A$. Then,  by definition $k\ge r(A)$.  Therefore, the second successor $A''$  of $A$ corresponds to the  first return of $f^k(c)$ to $A$, where $c$ is the 
critical point in $A$. Thus,  there exists  $k'$ such that $f^{k'}(f^k(c))\in A$ and  $f^{k+k'}(A'')=A$ (since $A''$ is successor of $A$, he cannot meet the critical point $c$ more than twice). Thus $k'\ge r(f^k(A''))\ge  r(A)$. By definition of $\sigma(A)$,  we conclude that $\sigma(A)\ge  k+k'$ and the result follows.
\cqfd

\subsection{Properties of the enhanced nest}

A  tool particularly well adapted for studying  the   persistently recurrent critical points is  {\it the enhanced nest} introduced by  Koslovski-Shen-vanStrien (see~\cite{KSS,QY,TY}). In our  situation we will use a slightly simplified version of the enhanced nest (see~\cite{TLP}). If the critical point $c$ is  persistently recurrent, the enhanced nest  consists in two sub-nests $(K_n),(K'_n)$ of  the critical nest 
$(P_j(c))$. The precise  construction  of these nests will be explained  in section~\ref{s:KSS}.
For the moment, we give  some of their properties.

\begin{remark}The construction in~\cite{KSS}, \cite{QY} or  \cite{TLP} has the feature  that the annuli between two consecutive pieces of the subnests are non-degenerate. However this property  is really not necessary for the construction itself. It is needed only afterwards  for the  control  of the moduli of these annuli. Here we do the construction without knowing  about the annuli. Then we prove that deep enough in the nest  the annuli   are non degenerate.
\end{remark}

The two sub-nests are constructed using  three operators called $\AA$, $\BB$ and $\DD$. As it is for $\DD$, $\AA$ and $\BB$ are some pull-back by  a conveniently  chosen  iterate of $f$. The operators  act on the set of all critical puzzle pieces and 
  "have bounded degree". The operators $\AA$ and $\BB$ are closely related and used for  producing  an annulus that avoids the closure of the orbit of the critical points in $perCrit(x)$. On the other hand, $\DD$  which is  just the last successor map (Notation~\ref{n:successor}) is used for giving  long iterates of bounded degree. 
  \vskip 0.5em
  {\it The subnests $K_n$, $K'_n$ are  defined inductively by $K_n=\AA\DD^{\tau}(K_{n-1})$ and $K'_n=\BB\DD^{\tau}(K_{n-1})$, starting from some $K_0$ containing the critical point of $perCrit(x)$ we focus on. Here $\tau $ is an arbitrary integer.}
    \vskip 0.5em

  \noindent
   In~\cite{TLP} we explain that taking $\tau=b+1$ is enough and that for $b\ge 2$, $\tau=b$ works also. 

\vskip 0.5emWe give here some properties of $\AA$ and $\BB$ that will be proved in section~\ref{s:KSS} (see also~\cite{KSS} or~\cite{TLP}).
Recall that $b$ is the number of critical ends and that $\delta$ is the maximal degree of $f$ over the critical ends. 
\begin{definition}Denote by $\displaystyle \PPc=\ol{\Cup_{c\in \omega Crit(c)}\Cup_{n\in \N}f^n(c)}$ generated by the critical points that are accumulated by $c\in perCrit(x)$. Note that  this set does not depend on the choice of the critical point in $perCrit(x)$.
\end{definition}

{\it Notice that in the post-critical orbit only the critical point is to avoid. However, if a critical point  belongs to  $c'$ the image of a puzzle piece $P_n(c)$,  necessarily 
 $c'\in \omega Crit(c)$. Hence, it is enough to control the presence of points of $\mathcal {P}_{\omega Crit(c)}$ in the  puzzle pieces. Also, any point of  $\mathcal {P}_{\omega Crit(c)}$ in a puzzle piece implies that an iterate of a critical point is inside.}

\begin{proposition}\label{p:pteAB} Let $I$ be a puzzle piece containing a critical point~$c$. Then $\AA(I)$ and $\BB(I)$ are  puzzle pieces with the following properties\,:

\begin{enumerate}
\item $c\in\AA(I)\subset\BB(I)\subset I$ and $\BB(I)\setminus \AA(I)$ is disjoint from the set $\PPc$\,;
\item there exist integers  $b(I),\  a(I)$ such that $f^{b(I)}(\BB(I))=I$ and $f^{a(I)}(\AA(I))=I$\,;
\item $\#\{0\le j< b(I)\mid c\in  f^j(\BB(I))\}\le b$, and   the degree of $ f^{b(I)}:\BB(I)\to I $ is less than $\delta^{b^2}$\,; 
\item $\#\{0\le j< a(I)\mid c\in  f^j(\AA(I))\}\le b+1$,  and   the degree of  $  f^{a(I)}:\AA(I)\to I $ is less than $ \delta^{b^2+b}$. 
\end{enumerate}
\end{proposition}
 \vskip 0.5em
\noindent  This Proposition has several  consequences on the  two sequences $(K_n),(K'_n)$. Let $h_n$ (resp. $h'_n$) denote the height of $K_n$ (resp. $K'_n$) in the nest\,: $K_n =P_{h_n}(c)$ and $K'_n =P_{h'_n}(c)$

\begin{corollary}\label{c:pteKSS}   The enhanced nests $(K_n),(K'_n)$ have the following properties\,:

\begin{enumerate}\item $K_n\subset K'_n\subset K_{n-1}$\,;
\item $K_n, K'_n$ are both  mapped to $K_{n-1}$ by some iterate of $f$\,: there exist integers    $p_n,p'_n$  such that\,: $f^{p_n}(K_n)=K_{n-1}$, $f^{p'_n}(K'_n)=K_{n-1}$\,;
\item $deg (f^{p_n}\from K_n\to K_{n-1})\le C$ and $deg (f^{p'_n}\from K'_n\to K_{n-1})\le C $ where $C$ depends only on $b$ and $\delta$\,;
\item $K'_n\setminus K_n$ is an annulus (possibly degenerate) that is disjoint from the  set $\PPc$\,;
\item $h(K'_n)-h(K_n)\ge r(K_{n-1})$ the return time in $K_{n-1}$. Moreover,   $r(K_{n+1})\ge 2^{\tau} r(K_n)$. So, $h(K'_n)-h(K_n)$ tends to infinity.
\end{enumerate}
\end{corollary}
\proof  Point 1) follows directly from Proposition~\ref{p:pteAB}.1 and the fact that  for any piece $J$, $\DD(J)\subset J$.  
 
Point  2)  follows directly from  Proposition~\ref{p:pteAB}.2.
Then Point 3) is exactly the second part of Proposition~\ref{p:pteAB}.1.
The point 4) follows from points 3. and 4. of  Proposition~\ref{p:pteAB} and  Remark~\ref{r:successor}. 
To  prove 5), note that $f^{p'_n}(c)$ and $f^{p_n}(c)$ are both in $K_{n-1}$, thus  
$p_n-p'_n\ge r(K_{n-1})$. The result follows from the fact that $p_n=h(K_{n-1})-h(K_n)$ and  $p'_n=h_{n-1}-h'_n$. The rest of the statement  follows from 
Corollary~\ref{c:2successors} since $K_{n+1}\subset \DD^\tau(K_n)$ so that $r(K_{n+1})\ge  r(\DD^\tau(K_n))$, that is bigger than $ 2^\tau r(K_n)$ .
\cqfd
In section~\ref{s:KSSproof} we give the proof of the following Lemma  that is a consequence of Lemma~\ref{l:2successors}\,:
\begin{lem}\label{l:doublement}
$p_n\ge 2p_{n-1}$.
\end{lem}
\begin{remark}\label{r:doublement}Lemma~\ref{l:doublement} is fundamental to understand the power of the construction. 
We get this way  long iterates of bounded degree\,:  let  $t_n$ be  the "time" necessary to reach  $K_0$, {\it i.e.}  such that $f^{t_n}(K_n)=K_0$. Then,  $t_n=p_n+\cdots+p_1\le  p_n+p_n/2+\cdots+p_n/2^{n-1}< 2p_n$, so that the  last step ($f^{p_n}\from K_n\to K_{n-1}$) 
takes more than half of the global time  and has a degree bounded by a constant  $C$   independent on $n$ (see point 4.).
Point 5 implies that the  annulus is large  in terms of the height.
\end{remark}
\subsection{Estimates on the moduli and bounds on the degrees}
Our next goal is to control the modulus $\mu_n$ of the annulus $K'_n\setminus K_n$. 
\begin{lem}\label{l:start}There exists two puzzle pieces $\tilde P, \tilde Q$ containing $c$ such that $\ol {\tilde Q}\subset \tilde P$. 
\end{lem} 
\proof Recall that $c\in \omega Crit(x)$ and that $\omega (x)$ contains no periodic point. 
Assume to get a contradiction that $P_0(c)\setminus P_n(c)$ is degenerate for every $n\ge 1$. Then $End(c)=\cap \ol P_n(c)$ contains an eventually periodic point $y\in \partial U$. By extension of notations, $End(c)=End(y)$. Corollary~\ref{c:kiwi} implies that $End(y)\cap \partial U=\{y\}$, so $y\in \omega(c)$ . This contradicts the fact that  $\omega (x)$ contains no periodic point.
\cqfd
\begin{definition} We define the`` {\it simplified enhanced nests}'' starting from $K_0=\tilde Q$.
\end{definition}
\begin{corollary}\label{c:Knondeg} For the  simplified enhanced nests, the annulus $K'_n\setminus K_n$ is non degenerate  for all large enough $n$.
\end{corollary}\proof
Since $f^{t_n}(K_n)=K_0$, pulling back by the good inverse branch of  $f^{t_n}$, we get a non degenerate annulus between a puzzle piece containing  $c$ called $R_n$ and $K_n$. The difference of height is constant\,: $h(K_n)-h(R_n)=h(Q)-h(P)$. Thus, since 
$h(K_n)-h(K'_n)$ tends to infinite, $h(K'_n)<h(R_n)$ for large $n$ so that 
$K'_n\supset R_n$. It follows that $K'_n\setminus K_n\supset R_n\setminus K_n$, so that it is non degenerate.\cqfd
We can find small pieces in $K_n$ on which  long iterates have bounded degre. Moreover, they are pull-back of $K_n$.
\begin{lem}\label{l:bdddegree}
Let $ z=f^\xi(c)\in K_n $ fro some $\xi$. Denote by $r_n$ the first entrance time of $f^{t_n}(z)$  in $K_n$. Let $A$  be the pull back of $K_n$  around $z$  by $f^{t_n+r_n}$.
Then, the degree of  $f^{t_n+r_n} \from A\to K_n $ is bounded by  $C_1:=2(C+\delta^b)$, where $C$ is the constant of Corollary~\ref{c:pteKSS}.
\end{lem}
\proof Since $z\in K_n$, the puzzle piece  $A:=P_{h_n+t_n+r_n}(z)$ is included in $K_n$. Therefore,   the degree of $f^{p_n}$ on $A$ is bounded by $C$ since the degree $f^{p_n}:K_n\to K_{n-1}$ is bounded by $C$.  
Now, take $s_n$ the first entry of $f^{p_n}(z)$ in $K_n$. By definition $f^{p_n}(z)$ enters 
in $K_n$ before $f^{t_n}(z)$ does, so $p_n+s_n\le t_n+r_n$.  Lemma~\ref{l:firstentrance}  implies  that $f^{s_n} \from  f^{p_n}(A)
 \to f^{p_n+s_n}(A) $ has  its degree  bounded by $\delta^b$. Now since $ f^{p_n+s_n}(z)\in K_n$ and since  the height of $ f^{p_n+s_n}(A)=P_{l_n}(f^{p_n+s_n}(z))$ is $l_n=h_n+t_n-p_n-s_n+r_n\ge h_n$,  the puzzle piece  $ f^{p_n+s_n}(A)\subset K_n$\,; therefore   $f^{p_n}$ has degree bounded by $C$ on it.
 Since $2p_n\ge t_n$ (Remark~\ref{r:doublement}) it follows that the  degree of $f^{t_n}$ on $A$ is bounded by $2C+\delta^b$.  Finally,  
the degree of $f^{r_n}: f^{t_n}(A)\to K_n$ is bounded by $\delta^b$ by Lemma~\ref{l:firstentrance}. \cqfd
\begin{corollary}\label{c:AB} Let  $A'$ be the pull back of $ K'_n$ around $z$  by $f^{t_n+r_n}$. Then $\mod(A'\setminus A)\ge  mod(K'_n\setminus K_n)/ C_1$ for $C_1=2(C+\delta^b)$.\end{corollary}
\proof Since $K'_n\setminus K_n$ does not intersect the set $\PPc$, the degree $f^{t_n+r_n}$ is the same on $A$ and on $A'$. The result then follows from the Lemma~\ref{l:bdddegree}. 
\cqfd

We explain now how we can compare the moduli between $K'_m\setminus K_m$ and another $K'_n\setminus K_n$.  It follows indirectly  from the fact that the height  $h(K'_m)-h(K_m)$ between $K'_m$ and $K_m$ is going to infinite, so  $K'_m\setminus K_m$ can contain sub-annuli which are related to    $K'_n\setminus K_n$. 

\begin{lem}\label{l:inclus} Let $\xi_n$ be the time such that  $f^{\xi_n} (K'_{n+2})=K_n$.
Then $f^{\xi_n}(K_{n+2})\subset A$ where $A$ is the pull back defined in Lemma~\ref{l:bdddegree}  with $z=f^{\xi_n}(c)$.
\end{lem}
\proof

Since  $f^{\xi_n}(K_{n+2})$ and $ A$  are both puzzle pieces containing $z=f^{\xi_n}(c)$, it suffices to compare $\alpha=\#\{0<j< t_n+r_n
\mid c\in f^j(A)\}$ with $\beta=\#\{0<j< p_{n+1}+
p_{n+2}-{\xi_n}
\mid c\in f^j(f^{\xi_n}(K_{n+2}))\}$.

The iterates $f^{i}(K_n)$ for $0\le i \le p_n$  meets at most $b+1+\tau$ times the point $c$ (by definition of $K_n=\AA\DD^\tau(K_{n-1})$ and using  Proposition~\ref{p:pteAB}). We can apply  Lemma~\ref{l:bdddegree} (and its proof), since   $f^{p_n}(z)$ returns in  $K_n$ because $c\in \omega Crit(c)$. Let $s_n$ be  the return time of 
$f^{p_n}(z)$ in $K_n$,  then  the images of $P_{h_n}(f^{p_n}(z))$
by $f^{i}$ for $0\le i < s_n$ never contain $c$. 
Therefore using that $p_n+s_n\le t_n+r_n$ we obtain  $\alpha\le 2(b+1+\tau)$.

On the other side, the number of iterates $p_{n+1}+
p_{n+2}-\xi$ to  bring $f^\xi(K_{n+2})$ to $K_n$  is less than $\sigma(K_n)(\beta-1)$
since $\sigma(K_n)$ (the number of  iterates for the last successor) is the largest time a point in $K_n$ take to come back to $K_n$.
It is also exactly the difference of height between $K_n$ and $f^\xi(K_{n+2})$ which is equal to the difference of height between $K'_{n+2}$ and $K_{n+2}$. Recall that 
$K'_{n+2}=\BB(\DD^\tau(K_{n+1}))$ and $K_{n+2}=\AA(\DD^\tau(K_{n+1}))$, so that this difference of height is bigger than $r(\DD^\tau(K_{n+1}))$, 
the return time in $\DD^\tau(K_{n+1})$. Moreover,  $r(\DD^\tau(K_{n+1}))\ge 2^\tau r(K_{n+1})$ ($\ge 2^{2\tau}r(K_{n}) $ by Corollary~\ref{c:pteKSS}). 
Now, since $K_{n+1}\subset \DD^\tau(K_{n})$,  we obtain that 
$r(K_{n+1})\ge r(\DD^\tau(K_{n}))\ge 2^{\tau-1} r(\DD(K_{n}))$. Then it follows from Remark~\ref{r:returntime} that    $r(\DD(K_{n}))\ge \sigma(K_n)$.
Hence $\beta\ge 2^{2\tau-1}$. 
We can summarise this discussion in\,:
$\sigma(K_n)(\beta-1)\ge p_{n+1}+p_{n+2}-\xi\ge h(K'_{n+2})-h(K_{n+2})\ge r(\DD^\tau(K_{n+1}))\ge 2^{2\tau_1}\sigma(K_n)$.

It is easy to see that for $\tau\ge b+1$ we obtain that $\beta\ge \alpha$, and the result follows.
\cqfd
\begin{corollary}\label{c:inclus}
Then $\mod(K'_{n+2}\setminus \ol K_{n+2})\ge \mod (K_n\setminus A)/C^2\ge  mod(K'_n\setminus \ol K_n)/ C'$ where $C$ and $C'$ are independant on $n$ and $\tau$.
\end{corollary}
\proof The map $f^\xi :   K'_{n+2}\to K_n$ has degree bounded by $C^2$ and the annulus $A'\setminus A$ is included in $K_n\setminus f^\xi(K_{n+2})$, so that 
$\mod (K'_{n+2}\setminus K_{n+2})\ge \mod (K_n\setminus f^\xi(K_{n+2}))/C^2\ge \mod (K_n\setminus A)/C^2\ge \mod (A'\setminus A)/C^2$. Finally, by Corollary~\ref{c:AB}, we get  
$C'=2(C+\delta^b)C^2$.
\cqfd
\subsection{ The Kahn-Lyubich Lemma}
The Covering Lemma due to Kahn and Lyubich  (see~\cite{KL}) is a very powerful tool\,:
it gives an ``estimate'' of the modulus of the pre-image of an annulus under a ramified covering  $g$ when one  has some control on the degree over some sub-annulus.
 We state  the Theorem now and refer   the reader to~\cite{KL}  for the proof.

\begin{theorem}[( The Kahn-Lyubich Covering Lemma)]
Let $g:U\to V$ be a degree $D$ holomorphic ramified covering. For any $\eta>0$ and $A,A',B,B'$ satisfying :
\begin{itemize}
\item
$A\subset \subset A'\subset \subset U$ and $B\subset\subset B'\subset \subset V$
are all disks\,;
\item  $g$ is a proper map from
$A$ to $B$, and from $A'$ to $B'$ with degree $d$\,;
\item  $mod(B'\setminus B)\ge \eta \ mod(U\setminus A)$.
\end{itemize}
 There exists 
$\epsilon=\epsilon(\eta,D)>0$ such that  $$mod(U\setminus A)>\epsilon\ 
\hbox{ or }\ mod(U\setminus A)>\frac \eta{2d^2}\  mod(V\setminus B).$$
\end{theorem}


\subsection{Proof of Theorem~\ref{th:fatoulc}}
 Using  the  results  obtained in previous sub-sections on the bound of the  degree of iterates of $f$ and the comparison on the moduli we get,   we now apply  Kahn-Lyubich  Covering Lemma  to prove.
 
 \begin{lem}\label{l:liminf}  The modulus $\mu_n$ of $K'_n\setminus K_n$ satisfies\,:  $\liminf \mu_n>0$.
\end{lem}
  
  \proof  Recall that  there exists some $N\ge 0$ such that $K'_n\setminus K_n$ is non degenerate for $n\ge N$ (Corollary~\ref{c:Knondeg}). To simplify the notations let us asume that $N=0$.
  Then the proof goes by contradiction. We assume  that  $\liminf \mu_n=0$. 
  Therefore  there exists a sequence $k_j\to \infty$ such that   $\mu_i\ge \mu_{k_j}$  for every $i\le k_j$   and $ \mu_{k_j}\to 0$.
    Since we already have indexes, we will fix some $n:=k_j-2$, so that $\mu_{n+2}\le \mu_i$ for $i\le n+2$. We would like to  apply the  Kahn-Lyubich Covering Lemma with $U=K_{n}$ and $V=K_0$, but then the  constant $\epsilon(\eta,D)$  would depend on the degree $D$ of the iterate $f^{t_n}$ from $K_n$ to $K_0$ and this degree goes to infinity (recall that  $f^{t_n}(K_n)=K_0$).
   
   Therefore,  we apply  the Covering Lemma    with $U=K_{n}$ and $V=K_{n-Z}$ for some integer $Z$. Let $\eta=1/(C^2C_1)$. We show that for 
    $Z>\frac{2C^2C_1^3}{\eta}$ the second case of the conclusion of the Covering Lemma cannot be satisfied.
 The set $A$, resp. $A'$, defined in Lemma~\ref{l:bdddegree}, resp. in  Corollary~\ref{c:AB},  is   the   pullback   of $K_n$, resp. $ K'_n$ around $z$  by $f^{t_n+r_n}$. 
Recall that $\xi$ is the time such that  $f^\xi (K'_{n+2})=K_n$,  that $z=f^\xi(c)$ and  that $r_n$ the first entrance time of $f^{t_n}(z)$  in $K_n$ ($r_n$ is well defined since $c\in \omega Crit(c)$). 
Let $z_n$ be the time such that $f^{z_n}(K_n)=K_{n-Z}$. Define $B:=  f^{z_n}(A)$ and  $B':=f^{z_n}(A')$. The puzzle pieces $B$, resp. $B'$   is  the   pullback   of $K_n$, resp. $ K'_n$ around $f^{z_n}(z)$  by $f^{b_n}$ where $b_n=t_n-z_n+r_n$.

We verify now the hypothesis of the Covering Lemma.
The map $g:=f^{z_n}\from U\to V $ is a covering of degree bounded by $ ZC$ (by Property 4. of Corollary~\ref{c:pteKSS}) and the puzzle pieces, $A\subset \subset A'\subset \subset U$ and $B\subset\subset B'\subset \subset V$
are all disks. Then $g$ is a proper map from
$A$ to $B$ and from $A'$ to $B'$\,; its    degree is  bounded by $ C_1$ since the degree of the maps $f^{t_n+r_n}: A\to K_n$ and $f^{t_n+r_n}: A'\to K'_n$ is bounded by $ C_1=2(C+\delta^b)$ independent of $n$  (Lemma~\ref{l:bdddegree}). Therefore, the maps $f^{b_n}\from B\to K_n$  and  $f^{b_n}\from B'\to K'_n$  have degree also bounded by $C_1$ and  $\mod(B'\setminus B)\ge \mod( K'_n\setminus  K_n)/C_1$.  
Since $\mu_n\ge \mu_{n+2}$ by assumption,  we see that   the  last hypothesis  is satisfied (using  Corollary~\ref{c:inclus})\,:
$$\mod(B'\setminus B)\ge \frac{\mu_{n+2} }{C_1}\ge\frac{ \mod(U\setminus A)}{C^2C_1}=\eta \mod(U\setminus A)\quad \hbox{ for } \eta:=\frac1{C'},  \hbox{ with } C'=C^2C_1  .$$

Thus the Kahn-Lyubich Lemma impies that  either    $\mod(U\setminus A)>\frac{\eta}{2d^2}\mod(V\setminus B)$ or there exists $\epsilon =\epsilon (\eta, D)$ independent on $n$ such that 
$\mod(U\setminus A)>\epsilon$.

We will prove that the first inequality  cannot arise  for   $Z>\frac{2C^2C_1^3}{\eta}$. The reason is that   the annulus $V\setminus B$ contains the pull back of $B$ around $f^{z_n}(z)$ of the annuli  $K'_{n-i}\setminus K_{n-i}$ for $0\le i\le Z$  by appropriate iterates of $f$.
Take the first time $f^{t_n}(z)$ enters in $ K_{n-i}$, let $B_i$ be the corresponding   pull back   around $f^{z_n}(z)$, and $B_i'$ the one of  $ K'_{n-i}$ . 
Applying Lemma~\ref{l:bdddegree} in this case (for the index $n-i$), we get   that 
 $\displaystyle \mod(B'_i\setminus B_i)\ge \frac{\mu_{n-i}}{C_1}$. 
 Hence, using that $\mu_{n+2}\le \mu_j$ for $j\le n+2$ and that   $d\le C_1$, we 
 obtain  $$\displaystyle \mod(V\setminus B)\ge \sum_{i=0}^Z \mod(B'_i\setminus B_i)\ge \frac{Z\mu_{n+2}}{C_1}.$$
 Using the Kahn-Lyubich Lemma, it follows that  
 $$\displaystyle  \mod(U\setminus A)>\frac{Z\eta}{2C_1^3} \mu_{n+2}.$$
 Combining this  inequality with the one   given in Corollary~\ref{c:inclus}\,:   $\displaystyle \mu_{n+2}\ge \frac{\mod(U\setminus A)}{C^2}$,   leads to $\displaystyle \mu_{n+2}>\frac{Z\eta}{2C^2C_1^3} \mu_{n+2}$, and by the choice of $Z$ one gets the contradiction that $\mu_{n+2}>\mu_{n+2}$.

Finally,  there exists $\epsilon=\epsilon(\eta,D)>0$ independent of $n$ such that $ \mod(U\setminus A)>\epsilon$. Hence, 
$$\mu_{n+2}\ge \mod(U\setminus A)/C^2>\epsilon/C^2.$$ Since this quantity is independent on $n$  we obtain that $\liminf \mu_i>0$ and get the contradiction. 
\cqfd

 \begin{corollary}\label{c:diam} The diameter of $K_n$ tends  to $0$.
\end{corollary}
\proof The annuli $A_n=K'_n\setminus K_n$ are disjoint and essential in the annulus  $A:=K_0\setminus (\cap K_n)$. By  Lemma~\ref{l:liminf}  there exists some $\epsilon'>0$  such that for large $N$, every $\mu_{n}\ge \epsilon'/2$ for $n\ge N$. Therefore  by Gr\"otzsch inequality (see~\cite{Ah}),  $\sum \mod(A_i)\le \mod(A)$ and in particular $\mu(A)=\infty$.
To conclude we use the following classical  characterization (see~\cite{Ah})\,: A continuum $K$ contained in a disk $D$ is reduced to a single  point if and only if $\mod(\D\setminus K)=\infty$.
\cqfd

\begin{corollary}\label{c:Ex} Let  $x$ be a point of $ \partial U$  satisfying  $\omega Crit(x)= perCrit(x)$. Then $End(x)=\{x\}.$
\end{corollary}
\proof Let $c_0$ be a critical point in $\omega Crit(x)$. By Lemma~\ref{l:start} there exist two puzzles pieces $\tilde Q$ and $\tilde P$ such that $c_0\in \tilde Q$ and $\ol{ \tilde Q}\subset \tilde P$. We construct the simplified enhanced nest around $c_0$ starting with $K_0=\tilde Q$. Then Lemma~\ref{l:liminf} implies that $\liminf \mu_n>0$ so $End(c_0)=\cap K_n$ reduces to $c_0$. 

For each $n\ge 0$, we consider the first entrance time of $x$ to $K_n$\,:  there exists $k_n$ such that $f^{k_n}(x)\in K_n$. Then the puzzle piece $P_{h'_n+k_n}(x)$ (resp. 
 $P_{h_n+k_n}(x)$) is mapped to $K'_n$ (resp. $K_n$) by $f^{k_n}$. These two maps 
 $ f^{k_n}_{\vert P_{h'_n+k_n}(x)}$ and  $ f^{k_n}_{\vert P_{h_n+k_n}(x)}$ have degree bounded by $\delta^b$.  Therefore the annulus $A_n(x):=P_{h'_n+k_n}(x)\setminus \ol{P_{h_n+k_n}(x)}$ has modulus $\mu_n(x) \ge \frac 1{\delta^b}\mu_n$. Moreover, the annuli 
 $A_n(x)$ are disjoint and nested around $x$ since $K'_{n+1}\subset  K_n$.
 Therefore we obtain that $mod(P_0(x)\setminus End(x))\ge \sum \mu_n(x)\ge  \frac 1{\delta^b}\sum \mu_n=\infty$, as in previous Corollary. The result follows.
\cqfd

\section{The operators $\AA$ and $\BB$}\label{s:KSS}
 This section is devoted to the definition of the operators $\AA$ and $\BB$.
 
 The 
definition    and the proofs are exactly the same as in~\cite{KSS, TLP,QY}. However, this construction is always presented with sequence of puzzle pieces that do not touch at the boundary. 
Since in our situation the puzzle pieces may touch, we will give here the details of the construction.

{\bf Important\,:} Through all this section  we will call ``annulus'' the difference $U\setminus U'$ between two open disks (or between one open disk and the closure of a smaller one) whose boundaries possibly touch at finitely many points.
Here these disks will be always puzzle pieces.  
(Notice that with this definition, an annulus is not always connected).

\subsection{Understanding the pullback to avoid the part of the post-critical set}
We fix some point $c_0\in perCrit(x)$ where $x$ is a given point of $\partial U$. We  look after  annuli around $c_0$ that avoid the post-critical set  $\PPo$.
By assumption~\ref{a:nocrit},
if  $y\in\PPo\cap P_0(c_0)$ then $y=f^j(c)$ for some $c\in \omega Crit(x)$ then $c\in \wc(c_0)$.  Therefore we focus on the  critical points of $\omega Crit(c_0)$. 

Before entering into the details of their definition, we sketch briefly the construction of the operators $\AA$ and $\BB$. These operators act on the set of puzzles pieces containing  a given critical point, here  we call it $c_0$.  Let $I$ be such a puzzle piece, consider the pullback of $I$ by the iterates $f^k$ and $f^{k'}$  ($k,k'>0$) corresponding to the first and to the second entrance of the orbit of a critical point $c\in\omega Crit(c_0)$. Doing this for each critical  point $c\in \omega Crit(c_0)$  by induction with a union of puzzle pieces instead of $I$, we obtain a finite set of pieces denoted by $P_c$. 
Then, we take  another   pullback inside $P_c$ called $P'_c$. Now, consider  among all the successors of $I$, those mapped to  $P_c$. Then,   $\BB(I)$ is defined as one of these successors satisfying   some maximality property. The piece  $\AA(I)$  is defined as  the pullback by   $f ^{b(I)} 
:\BB(I)\to I$, of $W$,  the first pullback of $I$ around  $f ^{b(I)}(c_0)$ ($W=\LL_{f ^{b(I)}(c_0)}(I)$).

When the map $f$ has  only one critical point $c_0$,  $\BB(I)$ is simply the last successor $\DD(I)$ and $\AA(I)$ the pullback in $\BB(I)$  (by $f^{\sigma(I)}$) of $W$ the first pullback of $I$ around $f^{\sigma(I)}(c_0)$.

We enter now into some important remarks  and properties about the pullback   by iterate corresponding to the first entrance in a puzzle piece.

\begin{definition}\label{d:Lz}\begin{itemize} \item If $P$ is a puzzle piece and $z\in \C$,  
 denote by $\LL_z(P)$ the puzzle piece containing $z$ that is mapped by $f^k$ to $P$,
 where $k>0$ is the first entrance time of $z$ in $P$, if it exists. 
 \item 
 If $H$ is a finite union of puzzle pieces,    let  $k>0$ be the first entrance time of $z$ in $H$ and let  $P$ the component of $H$ which contains $f^k(z)$.  Then define  $\LL_z(H):=\LL_z(P)$. 
 \item We call    $\LL_z(P)$, resp. $\LL_z(H)$,  {\it the first pull back} of $P$ (resp. of  $H$) around $z$. 
 \item
 Let $f^k(z)$   the second entrance of $z$  in $H$. We 
 note $\LL''_z(H)$ the pullback by $f^k$  around $z$  of the puzzle piece $P\subset H$  containing 
$f^k(z)$. We  call it  {\it the second pull back} of $H$ around~$z$. 
 \end{itemize} \end{definition}

\begin{remark}Note that  $\LL''_z(H)$ is also the pullback by 
 $f^l$  around $z$ of $\LL_{f^l(z)}(H)$, where $0<l<k$ label the successive  first and second entrance of $z$ in $H$. Thus   $\LL''_z(H)=\LL_z(\LL_{f^l(z)}(H))$.
 \end{remark}
 
\begin{lem} \label{l:pullbackdisjoints} Let $P$ be a puzzle piece. If $z'\notin \LL_z(P)$, then 
$\LL_z(P)\cap \LL_{z'}(P)=\emptyset$. 
\end{lem}\proof
We argue by contradiction. 
Let $k,k'$ be the first entrance time of $z,z'$ in $P$ so that $f^k(\LL_z(P))=P$ and $
f^{k'}( \LL_{z'}(P))=P$. By assumption,  $\LL_z(P)\cap \LL_{z'}(P)\neq \emptyset$, so $\LL_z(P)\subset \LL_{z'}(P)$ since they are puzzle pieces and since $z'\notin \LL_z(P)$. It implies that the levels of the puzzle pieces satisfy $h(\LL_z(P))> h(\LL_{z'}(P))$. Therefore, $k>k'$ since $k=h(\LL_z(P))-h(P)$ and $k'=h( \LL_{z'}(P))-h(P)$. The contradiction comes from the fact that  $k\le k'$ since it is the first entrance time in $P$ and $f^{k'}(z)\in P$.
\cqfd

The property of  Lemma~\ref{l:pullbackdisjoints}  extends  to finite unions  of puzzle pieces $H$ which   satisfy an extension property. In the litterature there exists already similar  properties of $H$ called {\it nice} or   {\it strictly nice}\,: recall that $H$ is  {\it nice}, resp.   {\it strictly nice}, if for every 
$n\ge 1$ and  every $z\in \partial H$, the iterate $f^n(z)\notin H$, resp. $f^n(z)\notin \ol H$.
The property we will use is of being {\it decent}. 

\begin{definition}Let $H$ be a  union  of finitely many open Jordan disks. One  says that $H$ is  {\it decent} if\,: \begin{enumerate}\item $H$ is nice\,;
\item no connected component of $H$ is mapped exactly to a component of $H$.
\end{enumerate}
\end{definition}
One can visualize this as `` the iterates of any component of $H$  is strictly bigger (not equal)   than a component of $H$ or doesn't intersect $H$''. Note also that we  do not require  to have an non degenerate annulus between an image of a component of $H$ and the component of $H$ contained in this image.

We are going to prove that a  finite union of disks which are decent satisfy the conclusion of  Lemma~\ref{l:pullbackdisjoints}. We begin with next Remark which follows from the definition of decent.  
\begin{remark}\label{r:weakinclus}If $H$ is a  finite union of puzzle pieces that is decent and  if $z\in H$, 
then  $\LL_z(H)\subset H$ for $z\in H$.
\end{remark}
\proof Suppose, in order  to get a contradiction, that    $\LL_z(H)\nsubseteq H$.
%
 Let $P$ be the connected component of $H$ which contains $f^k(z)$, the first return of $z$ in $H$, and  let $R$ be the component of $H$ containing $z$. 
Since puzzle pieces are either disjoint or  nested, $R$ is contained in  $\LL_z(H)$. 
Therefore,  $f^k(R)\subset f^k(\LL_z(H))=P$. Since $H$ is nice, $f^k(R)=P$
and this contradicts the fact that $H$ is decent.
\cqfd

\begin{lem}\label{l:nice} If $H$ is a  finite union of puzzle pieces that is decent, 
then for any two points $z,z'\in \C$, the pieces    $\LL_z(H)$ and $\LL_{z'}(H)$  either are disjoint or equal. 
\end{lem}
\proof Suppose that   $\LL_z(H)$ and $\LL_{z'}(H)$   are not disjoint.
Assume, in order to get a contradiction,  that $\LL_z(H)\neq \LL_{z'}(H)$, then  say  $\LL_z(H)\supset\LL_{z'}(H)$ since both are puzzle pieces. Let $k,k'$ be the respective entrance time of $z,z'$  in $H$. Then $f^k(z')\in f^k(\LL_z(H))\subset H$,  so  $k'\le k$.
Now $f^{k'}(\LL_{z'}(H)$ is equal to a connected component $P_i$ of $H$ so that $f^{k-k'}(P_i)\subset  f^{k}(\LL_z(H))$ which is  a connected component $P_j$  of $H$. Since $z\notin  \LL_{z'}(H)$ then   $k- k'\neq0$. As noticed during the proof of the Remark~\ref{r:weakinclus}, it is not possible that  $f^{k-k'}(P_i)\subset P_j$,  since  $H$ is decent. 
From this contradiction we get that  either   $\LL_z(H)=\LL_{z'}(H)$ or 
 $\LL_z(H)$ and $\LL_{z'}(H)$ are disjoint. 
 \cqfd
 
 \subsection{Construction of preferred puzzle pieces $(P_c)$ around the critical points} Let $x\in \C$ be such that $perCrit(x)\neq \emptyset.$
 Recall that $b$ is the number of distinct critical ends and $\delta $ is the maximum of the degree of $f$ over the critical ends  (see Notation~\ref{n:numbercritend}).  \begin{proposition} \label{p:Pc} Let $I$ be a puzzle piece around $c_0\in perCrit(x)$. There exists  two  sets   of puzzle pieces\,: 
 $$\Lambda:=\{P_c\mid  c\in \omega Crit(c_0)\}\hbox{ and } \Lambda':=  \{P'_c\mid  c\in \omega Crit(c_0)\}$$ with the following properties\,:
 \begin{enumerate}
 \item $\forall c\in \omega Crit(c_0)$,  $ c\in P'_c$ and $P'_c\subset P_c$\,;
 \item if $b=1$ then $P_c=I$ and $P'_c=\LL_{c_0}(I)$\,;
 \item else $P_c$ is a pullback of $I$ by some iterate $f^p$ and 
 $\left\{\begin{array}{l}\deg(f^p \from P_c\to I)\le \delta^{b^2-b};\vspace{0.12cm}\\
 \# \{0\le i<p\mid c_0\in f^i(P_c)\}\le b-1\ .
             \end{array}\right.$
\item The piece $P'_c$ is also a pullback of $I$. Moreover, if there is a point $z\in \PPo
 \cap (P_c\setminus  \overline {P'_c})$, then  there exists a   puzzle piece $V$ that   is included in  $ P_c\setminus  \overline {P'_c}$  and $f^k\from V\to P_{c'}$  is an homeomorphism.
\end{enumerate}
\end{proposition}

\noindent Recall that $\PPo$ is the closure of the forward orbit of the points in $\omega Crit(c_0)$.  Remark that if we note  $H:=\Cup_{c\in \omega Crit(c_0)} P_c$,
the piece $V$ given by  point 4)  is mainly $\LL_z(H)$ or some pullback of such a piece and the map $f^k\from V\to P_{c'}$ is   the first return  map  in $H$. 

\proof
We define the set by  induction on the number of critical points. Let  $H_0:=I$ and $J_0:=\LL_{c_0}(I)$.
\vskip 0.5 em
\noindent 1. {\it Assume that every  point of $\omega Crit(c_0)\setminus\{c_0\} $  enters in $J_0$ when   it enters in $H_0$ for  the first time}\,:

\noindent
 Then 
we take $P_{c_0}:=H_0$ and $P'_{c_0}:=J_0$,  for $c\in \omega Crit(c_0)$ we take   $P_c:=\LL_c(H_0)$ and $P'_c:=\LL_c(J_0)$. 
The property on the degree is clear for $P_c$ and for $P'_c$ (see Lemma~\ref{l:firstentrance}).

Assume that $b=1$, {\it i.e.} there is only $c_0$ in $Crit$. If $\PPo\cap (H_0\setminus J_0)\neq \emptyset$ then it should be   some iterate $f^r(c_0)$. Let  $V:=\LL_{f^r(c_0)}(H)$ where $H:=\Cup_{c\in \omega Crit(c_0)} P_c$, then
 $V\cap J_0=\emptyset$ (by Lemma~\ref{l:pullbackdisjoints}). Moreover,  $f^k\from V\to H$, the first return  map  in $H$,  is an homeomorphism because 
 $f^k$ has no critical points in $V$. 
 This step proves the  point 2)  of the Proposition.
\vskip 0.5 em
\noindent 2. {\it In the case $b>1$, we suppose that there exists $r\ge 1$ and two finite unions of puzzle pieces $J_r$, $H_r$ such that\,:
\begin{enumerate}
\item 
$ \{c_0,\cdots, c_r\}\subset J_r\subset H_r$\,; 
\item
$H_r$ is decent, the components of $J_r$ are of the form $\LL_c(H_r)$\,;
\item
every  $c\in \CC:=\omega Crit(c_0) 
\setminus \{c_0,\cdots, c_r\} $  enters  in $J_r$ when it enters   in $H_r$ for the first time.
\end{enumerate}}

\noindent Then we  define  $P_c$ and $ P'_c$ as follows:  for $0\le i\le r$, $P_{c_i}:=H_r^i$, resp.   $P'_{c_i}:=J_r^i$, is   the connected  component of $H_r$, resp. of $J_r$, containing $c_i$ and  for $r<i$  $P_{c_i}:=\LL_{c_i}(H_r)$,   $P'_{c_i}:=\LL_{c_i}(J_r)$. 

The proof of  point 4) of the Proposition goes as follows. Let  $z\in \PPo \cap (H_r\setminus  J_r)$ and let $V=\LL_z(H)$. If    $V=\LL_z(H^i_r)$ for some $i\le r$,  it is clear by Lemma~\ref{l:nice} that $V\cap J_r=\emptyset$ since the components of $J_r$ are of the form
$\LL_c(H_r)$ and $z\notin J_r$. 
Else,  $V=\LL_z(P_{c_i})$ for some  $i>r$. Suppose to get a contradiction that   $V\cap J_r\neq \emptyset$.  Then $V\supset L_{c_i}(H_r)$ (but they are not equal) for some $0\le i\le r$. Let $l,k$ be the respective entrance time of $z$ to $H$\,: $f^l(V)=P_c$ and $f^k(P_c)=H^j_r$ for some $j\le r$. Then $f^{k+l}(\LL_{c_i}(H_r))\subset H^j_r$ but they are not equal.
This contradicts the fact that $\LL_{c_i}(H_r)$ is the pullback of the puzzle piece of the first entrance to $H_r$ (so cannot be included in).
Note that since $V\cap J_r=\emptyset$ there is no critical points in $V$. Then, the first iterate $m$  such that $f^m(V)$ contains a critical point $c_i$ (by the recurrence assumption it is necessarily in $\omega Crit(c_0)$), is contained in the corresponding  connected component of $H$  (it cannot be bigger by Lemma~\ref{l:nice} since it is eventually mapped to such a component). Therefore, $f^m(V)$ is exactly the  component $P_{c_i}$ of $H$, so that $f^m\from V \to P_{c_i}$ is an homeomorphism.

Now we consider a point $z\in \PPo \cap (P_c\setminus  \overline {P'_c})$ for some $c\in \CC$.
Let $y$ be the first entrance of $z$ in $H_r$ and $V$ be the pullback of $\LL_y(H)$ around $z$.
 Note that $y\notin J_r$. 
 Indeed, the map from $P_c\setminus P'_c$ to $\H_r^j\setminus J_r^j$ is a non ramified covering since
 all the critical points are mapped into $J_r$. Hence, by the previous discussion, $\LL_y(H)$ is disjoint form $J_r$.
 And pulling back by the covering, we get that  $V\cap P'_c=\emptyset$.
 Moreover, the pullback by a non ramified covering of a disk is an homeomorphism, so $f^k\from V \to P_{c_i}$
 is an homeomorphism (using previous study).
 
 This achieves the proof of point 4).
We will see later that Point 3) is satisfied also. This will follow from the construction of $H_r$ and $J_r$.

\vskip 0.5 em
\noindent 3. {\it The construction of $H_r$ and $J_r$\,:}

 The  union of puzzle pieces $H_r,J_r$ with the properties stated in 2) are constructed by induction. Assume  that we have constructed a set $H_m$ union of  puzzle pieces $H_m^i$ around $c_i\in \omega Crit(c_0)$  for $0\le i\le m$ such that   $H_m$ is 
decent.
Let $J_m=\cup_{i=1}^m \LL_{c_i}(H_m)\subset H_m$.

If every critical point of $\omega Crit(c_0)\setminus \{c_0,\cdots, c_m\}$ enters in $J_m$ the first time it enters in $H_m$, then $H_m$ and $J_m$ have the properties required in 2) and we are done.  

\noindent
If not,  there is  a  point $c_{m+1}\in \omega Crit(c_0) 
\setminus \{c_0,\cdots, c_m\} $  that enters in $H_m\setminus J_m$ the first time it enters $H_m$\;, we
  set $H_{m+1}^{m+1}:=\LL''_{c_{m+1}}(H_m)$ and  $H_{m+1}:= J_m\cup H_{m+1}^{m+1}$. 
  Define $J_{m+1}:=\cup_{i=0}^{m+1} \LL_{c_i}(H_{m+1})$.

If we  prove that  $H_{m+1}$ is decent, it will follow by Remark~\ref{r:weakinclus}  that 
 $J_{m+1}\subset H_{m+1}$. After finitely many steps we will achieve the induction, the sets $H_m$ and $J_m$ will  have the properties required in 2).

 By definition, it is clear that $J_m$ is also  decent. 
Hence,   we only need to verify that no puzzle piece of $J_m$ can be mapped into 
 $H_{m+1}^{m+1}$. Let $k$ be the iterate such that $f^k(H_{m+1}^{m+1})$ is a $H_m$ puzzle piece.
 If $f^l(J_m)\subset H_{m+1}^{m+1}$ but does not coincide, then $f^{k+l}(J_m)$ is included in a $H_m$ puzzle piece but does not coincide with this piece. This is in  contradiction with the definition of $J_m$ that is the pullback of the $H_m$ puzzle piece of the first entrance.

Now consider  the puzzle piece $H_{m+1}^{m+1}$. Indeed, none of its  iterates 
can be contained in $H_{m+1}^{m+1}$ for level of piece reason.
Assume now, to get a  contradiction,   that for some $j$,  
$f^j(H_{m+1}^{m+1})\subset J_m$. Recall that $f^i(c_{m+1})\notin H_m$ for $0\le i\le k$ and $i\neq l$
where $l,k$ are respectively the first and the second entrance to $H_m$. Recall also that $H_{m+1}^{m+1}=\LL_{f^{l}(c_{m+1})}(H_m)$, so that it is disjoint from $J_m$ by Lemma~\ref{l:nice}. Moreover,  when it meets
$H_m$ for the second time, it is exactly along a  connected component of $H_m$, so it contains a  component of $J_m$. This contradicts the fact that some iterate is included in a $J_m$ component. 

\vskip 0.5 em
\noindent 4. {\it Degree properties :}
Let us prove now the second point of the Proposition.  We stop the induction at some step $r$.
For all $m\le  r$, any puzzle piece $H_m^i$ with $i\neq m$ is mapped back to a $H_{m-1}$ piece, and meets at most 
one time 
the  critical points of $\omega Crit(c_0)\setminus\{c_0,\cdots,c_{m-1}\}$ and only one in the set $\{c_0,\cdots,c_{m-1}\}$ so the degree is less than or equal to $\delta^{b-m+1}.$
The puzzle piece $H^m_m$ does twice  the turn, therefore the degree is $\le \delta^{2(b-m)}$.
Therefore the degree from $P_c$ to $I$ is less than the product of the $ \delta^{2(b-m)}$,  for $0\le m\le r$, which is $\le \delta^{b^2-b}$ since $r\le b-1$.

Now since at each step one meets $c_0$ at most once, the iterates from  $P_c$ to $I$ meet $c_0$ at most 
$r\le b-1$ times.
This finishes the proof of the Proposition~\ref{p:Pc}.
\cqfd

\subsection{Construction of $\AA$ and $\BB$}

In this section we give the precise definition of the operators $\AA$ and $\BB$, using the previous construction of the puzzle pieces $P_c,P'_c$. 
Let $I$ be a puzzle piece around $c_0$. To get an annuli avoiding $\PPo$ we take pullbacks and  ask to go through a puzzle piece of the set $\{P_c\mid c\in \omega Crit(c_0)\}$. So we want to consider the successors of $I$ (since we need annuli around $c_0$)
that goes through such puzzle pieces. But this set might be empty, as we have seen that we might meet several times the same critical point from $P_c$ to $I$.   We want to look at puzzle pieces $S$ containing $c_0$, that are mapped to $P_c$ and which meet 
the critical points at most twice before $P_c$. This corresponds to the notion of successor of $P_c$ generalised  because  we go back to $c_0$ and not to $c$. We  can also  use the notion of child as follows. Consider the  following collection of pieces (see Definition~\ref{d:Lz})
$$ \tilde \Lambda_{c_0}:=\{\LL_{c_0}(Q)\mid Q  \hbox{ is a child of  } P_{c} \hbox{ and } c\in \omega Crit(c_0) \}.$$ 

We detail  now the maximality property that characterizes $\BB(I)$ in  $ \tilde \Lambda_{c_0}$. Denote by 
$\tau_c(Q)$   the iterate such that $f^{\tau_c(Q)}(Q)=P_c$  when $Q$ is a child of $P_c$.
\begin{remark}  The set $\aleph:=\{\tau_c(Q) \mid Q\hbox{ is a child of }P_c ,\ c\in \omega Crit(c_0)\}$ is finite. 
\end{remark}
\proof It follows from the fact that  $c_0\in perCrit(x)$. Indeed, this  implies  that 
any piece $P_c$ with $c\in \omega Crit(c_0)$ has finitely many successors, and therefore finitely many children.\cqfd

\begin{definition}
Let $\tau=max \aleph$ and $c_\aleph$ a point of $\omega Crit(c_0)$ where the maximum  is reached and $Q_\aleph$ the corresponding child of $P_{c_\aleph}$\,:    
 $f^{\tau_\aleph}(Q_\aleph)=P_{c_\aleph}$.
Denote by $\wt c$ the critical point  of  $Q_\aleph$. \end{definition}

\begin{lem} By the definition of $\Lambda$ and $\Lambda'$, we get the property that  $f^\tau(\wt c)\in P'_{c_\aleph} $.  Moreover, if $Q_{c_\aleph}'$ is  the pullback containing $\wt c$ of $P'_{c_\aleph}$  by $f^\tau: Q_{c_\aleph}\to P_{c_\aleph}$, then $(Q_{c_\aleph}\setminus Q_{c_\aleph}')\cap \PPo=\emptyset$.
\end{lem}
\proof 
Assume in order to get a contradiction that  $f^\tau(\wt c)\notin  P'_{c_\aleph}$. 
 Then by the third point of Proposition~\ref{p:Pc}, there exists a puzzle piece 
$V$ containing  $f^\tau(\wt c)$ and an iterate of $f$ that is a homeomorphism from $V$ to some  piece of $\Lambda$, say $P_{\overline c}$ with $\overline c\in\omega  Crit(c_0)$.
Hence  the pullback $Q_0$ of $V$ by $f^\tau$  in $Q_{c_\aleph}$ around $\wt c$  is a child of this  $P_{\overline c}$\,:    $f^{\tau-1}$   from   $f(Q_{c_\aleph})$ to $P_{c_\aleph}$ is a homeomorphism and the map from $V\subset P_{c_\aleph}$ to $P_{\overline c}$ also. Note that the time from $Q_0$ to $P_{\overline c}$ is strictly larger  than $\tau$.  This contradicts the fact that   $\tau_{\overline c} \le \tau$.

Hence,  since  $f^\tau(\wt c)\in P'_{c_\aleph} $, we can define $Q_{c_\aleph}'$ the pullback   of $P'_{c_\aleph}$ containing $\wt c$. The  map    $f^\tau$ is a non ramified covering from $Q_{c_\aleph}\setminus Q_{c_\aleph}'$ to $P_{c_\aleph} \setminus P'_{c_\aleph}$, since $Q_{c_\aleph}$ is a child of  $P_{c_\aleph} $ and since there is  at most one critical point in each piece (here it is in $Q_{c_\aleph}'$). 

 Suppose, in order to get a contradiction, that  $\PPo\cap (Q_{c_\aleph}\setminus Q_{c_\aleph}')\neq \emptyset$. Let $z$ then be  a point  in $ \PPo\cap (Q_{c_\aleph}\setminus Q_{c_\aleph}')$,  image of a critical point $c$\,: $z=f^k(c)$. Since  $f^\tau(z)\in \PPo\cap( P_{c_\aleph}\setminus  P'_{c_\aleph})$,  there exists a puzzle piece $V$ containing $z$ and a point $\overline c\in \omega Crit(c_0)$ such that  the iterate from  $V$ to $P_{\overline c}$ is a homeomorphism. As before, the pullback of $V$ by  $f^\tau$ is a disk $V'$ in $Q_{c_\aleph}\setminus Q_{c_\aleph}'$ on which $f^\tau $ is an homeomorphism. Then, 
taking the pullbacks along the orbit  of $c$, the last iterate $f^i(\LL_{c}(V'))$ for $i<k$  which contains a critical point 
gives a child of   $P_{\overline c}$. The contradiction comes again  from the fact that  $\tau_{\overline c}$ is strictly greater than $\tau$.
  \cqfd

\begin{definition} Let $\BB(I):=\LL_{c_0}(Q_{c_\aleph})$.   Let $\AA(I)$ be the pullback by  $f ^{b(I)} 
:\BB(I)\to I$ of $W=\LL_{f ^{b(I)}(c_0)}(I)$.
\end{definition}

\begin{lem}\label{l:ABPf}  By construction,   $(\BB(I)\setminus  \AA(I))\cap \PPo=\emptyset$.
\end{lem}
 \proof
 There is an integer $n$ such that  $f^n(\BB(I))=Q_{c_\aleph}$. We will  prove that $f^n(\AA(I))\supset Q_{c_\aleph}'$\,; it will follow that any point of $(\BB(I)\setminus  \AA(I))\cap \PPo$ will have its image under  $f^n$  which is contained in $\PPo\cap (Q_{c_\aleph}\setminus Q_{c_\aleph}')$. But this set is empty   by the previous Lemma.  
 In order to prove that $f^n(\AA(I))\supset Q_{c_\aleph}'$,  notice that since $\PPo\cap (Q_{c_\aleph}\setminus Q_{c_\aleph}')=\emptyset$, the image $f^n(c_0)$ cannot be in $Q_{c_\aleph}\setminus Q_{c_\aleph}'$. So  since  this point belongs to both   $f^n(\AA(I))$ and $ Q_{c_\aleph}'$, these two pieces are nested.  Assume that  $f^n(\AA(I))\subset Q_{c_\aleph}'$. 
Then, $W\supset f^{\tau-n}(Q_{c_\aleph}')\subset I$. By construction,  $Q_{c_\aleph}'$ is a pullback of $I$. We get a contradiction,  from the fact that  $ f^{\tau-n}(Q_{c_\aleph}')$  will be mapped to $I$ before $W=\LL_{f ^{b(I)}(c_0)}(I)$. The Lemma the follows. 
 \cqfd

\subsection{ Proof of Proposition~\ref{p:pteAB}}
Recall the Statement of the Proposition\,:

\begin{proi}[{\bf\ref{p:pteAB}}] Let $I$ be a puzzle piece containing a critical point~$c$. Then the following holds\,:

\begin{enumerate}
\item $c\in\AA(I)\subset\BB(I)\subset I$ and $\BB(I)\setminus \AA(I)$ avoids the postcritical set $\PPo$\,;
\item There exists $b(I),\  a(I)$ such that $f^{b(I)}(\BB(I))=I$ and $f^{a(I)}(\AA(I))=I$\,;
\item $\#\{0\le j< b(I)\mid c\in  f^j(\BB(I))\}\le b$, and  $deg(  f^{b(I)}:\BB(I)\to I )\le \delta^{b^2}$\,; 
\item $\#\{0\le j< a(I)\mid c\in  f^j(\AA(I))\}\le b+1$,  and  $deg(  f^{a(I)}:\AA(I)\to I )\le \delta^{b^2+b}$. 
\end{enumerate}
\end{proi}
\proof
\noindent Point 1) of the Proposition is just Lemma~\ref{l:ABPf}.
Then point 2) follows  just from the  definition.
To prove  the degree properties of point 3)  and point 4),  recall that  the iterate of $f$ that maps $\BB(I)$ to $ P_{c_\aleph}$ has degree bounded by 
$\delta^b$ (Lemma~\ref{l:firstentrance}). Then, point 3) follows  from the fact that  the map from  $ P_{c_\aleph}$ to $I$
has degree bounded by $\delta^{b^2-b}$.  By Lemma~\ref{l:firstentrance}, the degree from $W$ to $I$ is bounded by $\delta$.

Finally, if $l $ denotes the time such that  $f^l(\BB(I))=P_{c_\aleph}$,    the pieces   $f^i(\BB(I))$ for $0\le i \le l$ meets $c_0$ only 
once since   $\BB(I)=\LL_{c_0}(Q_{c_\aleph})$ is  the first pull back around $c_0$ of $ Q_{c_\aleph}$ which is a child of $P_{c_\aleph}$ (Definition~\ref{d:Lz}), and   it happens only   in $\BB(I)$. Therefore from the construction of $P_c$ we get easily Point 3) and Point 4).

\section{The properties of the enhanced nest}\label{s:KSSproof}
This section is devoted to the proof of Lemmas~\ref{l:2successors} and~\ref{l:doublement} which were crucial in the proof of Lemma~\ref{l:liminf} via the Lemmas~\ref{l:bdddegree},~\ref{l:inclus} and its power is unlighted in Remark~\ref{r:doublement}.
\subsection{
Proof of Lemma~\ref{l:2successors} (see also~\cite{TLP})}
This result is similar to the one for cubic and  quadratic maps  obtained by Branner-Hubbard and Yoccoz (see~\cite{BH},~\cite{Hu} and~\cite{Mlc}).  It is a crucial step for the estimates in Lemma~\ref{l:doublement}  for instance.

\begin{lem}[(Previously called Lemma~\ref{l:2successors})]
 Let $c$ be any critical point of $f$.
 \begin{enumerate}
 \item
  If $End(c)$ is not periodic, then each puzzle piece $P_n(c)$ has at least two successors. 
 \item If $c\in \omega Crit(x) $  has an end  $End(c)$ which is  periodic and $x$ is a point of $\partial U$, then $\omega(x)$ contains a periodic point.
 \end{enumerate}
\end{lem}
\proof
Let $P$ be a puzzle piece containing $c$. Denote by $n_0$ its height\,: $P=P_{n_0}(c)$.  
Let  $k_0=0<k_1<\cdots<k_n<\cdots$ denote  the successive entrance time  of (the orbit of) $c$ in $P$.  
Note that  this  sequence is infinite since $c\in \omega Crit(c)$.  The first entrance time of $c$ to $P$ is  $k_1$,  denote by $Q$  the pull back of $P$ by $f^{k_1}$, it  is the first successor of $P$. Note that $Q=P_{n_0+k_1}(c)$. 
\begin{claima}\label{c:1} Either $P$ has  a second successor or  $k_i=ik_1$ for all $i\ge 0$ and  $f^{ik_1}(c)\in Q$.
\end{claima}
\proof Assume first  that there exists $k_i$ such that $f^{k_i}(c)\notin Q$. It means that $c$ which is in $P_{n_0}(f^{k_i}(c))$ is not in  $P_{n_0+k_1}(f^{k_i}(c))$. This implies that there is no  critical points  in  $P_{n_0+k_1}(f^{k_i}(c))$ (else there would be two critical points in $P_{n_0}(f^{k_i}(c))$ which contradicts assumption~\ref{a:3}). Then consider $Q'$ the    pull back of  $P$ containing $f^{k_i}(c)$ by $f^{k_{i+1}-k_i}\from P\to  f^{k_{i+1}-k_i}(P)$. If ${k_{i+1}-k_i}\ge k_1$, we just show that $Q'$ contains no critical points.   If ${k_{i+1}-k_i}< k_1$, if $Q'$ contains a critical point it should be $c$ (since $Q'\subset P$) but then $Q'$ is the first pull back of $P$ in $P$ which is impossible.
Therefore we can continue to pull back $Q'$ by the iterates of $f$ 
along the orbit $\{c,f(c),\ldots,f^{k_{i}}(c)\}$ until we reach a critical puzzle piece $R$ of depth  greater than that of $Q'$  that is greater than $n_0+k_1$.
If $c\in  R$  then $R$ is a second successor of $P$.
Else, $\LL_c(R)$ will be a second successor of $P$.

Now assume that  for every $i\ge 1$, $f^{k_i}(c)\in Q$. 
In particular, 
$f^{k_1}(c)\in Q$ and therefore   $f^{2 k_1}(c)\in f^{k_1}(Q)=P$. It follows that 
$k_2\le 2k_1$.  On the other side, since $f^{k_1}(c)\in Q$, it follows that $f^{k_2-k_1}(Q)=P$. In particular,  $f^{k_2-k_1}(c)\in P$, so $k_1\le k_2-k_1$.    Hence, 
$k_2=2k_1$ and by induction $k_i=ik_1$ for all $i\ge 0$.
\begin{claima}\label{c:2}Assume that $P$ has only one successor, that $k_i=ik_1$  for all $i\ge 0$ and  that $f^{ik_1}(c)\in Q$.
Denote by  $R$ the puzzle piece $P_{n_0+2k_1}(c)$, it is  the pull back by $f^{k_1}$ of $Q$ around $c$.  Then, for all $i\ge 0$,   $f^{ik_1}(c)\in R$.
\end{claima}
\proof 
The idea of the proof is similar as before. We consider the critical points that appear in the pieces $\{f^j(R)\mid 0\le j\le k_1\}$. We denote by $j_0=0<j_1<\cdots<j_{r+1}=k_1$ the iterates of $R$ containing a critical point that we  denoted by $c_l$ with $0<l\le r$, denote by $S_n:=f^{j_n}(R)$.

First we shall prove by contradiction that there is no critical point in the puzzle pieces 
$f^j(Q)$ for $j_r<j<j_{r+1}$.
Assume that a critical point $\tilde c$ belongs to $f^j(Q)$ with $j_r<j<j_{r+1}$. Let $k=j_{r+1}-j$. The point $f^k(\tilde c)$ is included  in $P\setminus Q$ and $\LL_{f^k(\tilde c)}(P)$ as well. Since the map $f^k : \LL_{\tilde c}\to P$ has only a critical point at $\tilde c$, 
if $\tilde S$ is the connected component of $f^{-k}( \LL_{f^k(\tilde c)}(P))$, the piece 
$\LL_c(\tilde S)$ is a second successor of $P$. Indeed, the each  map of the form  $f^l: \LL_x(S)\to S$ meets each critical point at most once. 

The puzzle  piece $S_r$ contains $c_r$ and is of height $n_r:=n_0+2k_1-j_r$.  The piece $S_r$ is the only pre-image of $Q$ in $f^{j_r}(Q)$  by $f^{k_1-j_r} \from f^{j_r}(Q)\to P$. 
Since $f^{k_{i+1}-(k_i+j_r)}$ maps $f^{j_r}(Q)$ to $P$ and $f^{k_i}(c)\in Q$ so that $f^{k_i+j_r}(c)\in f^{j_r}(Q)$, we deduce that $f^{k_i+j_r}(c)\in S_r$ for all $i\ge 0$ since $f^{k_{i+1}}(c)\in Q$.
With the same argument, we prove successively that   $f^{k_i+j_l}(c)\in S_l$, for $l$ from $r $ to $0$, which corresponds to  $f^{ik_1}(c)\in R$ for $i\ge 0$ (for $l=0$).
\begin{claima}\label{c:3}
Either $P$ has a second successor or $End(c)$ is periodic.  \end{claima}
\proof We have proved in Claim~\ref{c:1} that either $P$ has a second successor or the assumption of Claim~\ref{c:2} holds. 
Then after  iterating the same argument  on   $P_{n_0+jk_1}(c)$  for $j\ge 2$ as in  Claim~\ref{c:2}, we obtain that the critical point $c$ belongs to every  piece $P_n(f^{ik_1}(c))$. Hence, the 
entire nest $P_n(c)$ with $n\ge n_0+k_1$ is mapped  by $f^{k_1}$ to the same nest $P_n(c)$ with $n\ge n_0$.
\begin{claima}\label{c:4}
Assume $End(c)$ is periodic for $c\in\wc(x)$ with $x\in\partial U$ then  $x$  combinatorially accumulates a periodic point.  \end{claima}
\proof 
Since $c\in \omega Crit(x)$ with $x\in\partial U$, every puzzle piece of the nest $(P_n(c))$ intersects $U$. Consider the internal rays in $U$ of $\partial P_{nk_1}(c)$\,: they are of the form $R_U(\zeta_n),R_U(\zeta'_n)$, where $\zeta_n$ and $\zeta_n'$ are two  adjacent sequences converging to some angle $\zeta$. Since the nest is ``fixed '' by $f^{k_1}$, it follows that $R_U(\zeta)$ is also fixed by $f^{k_1}$\,;  so it converges to a point $y$ in $End(c)$ that is  fixed by $f^{k_1}$.    Hence the end $End(c)=\cap_n\ol{P_n(c)}$ contains a $k_1$-periodic point. Therefore $End(y)=End(c)$ and since $c\in\omega Crit(x)$,  $x$ accumulates combinatorially the point $y$. Therefore $y\in \partial U$ (because $x\in \partial U$) and since $End(y)\cap \partial U=\{y\}$ it follows that $x$ accumulates $y$ ({\it i.e.} $y\in \omega(x)$).
\cqfd

\subsection {
Proof of Lemma~\ref{l:doublement}}

Set   $I_n=\DD^\tau(K_{n-1})$, where $\DD(J) $ denote the last successor of a puzzle piece $J$, and $\tau$ is some number that can be for instance $b+1$. Recall that $K_n=\AA(I_n)$ and that $p_n, a(I_n), \sigma_\tau(K_n)$ are  the times to go respectively from $K_n$ to $K_{n-1}$, from $\AA(I_n)$ to $I_n$
 and from $I_{n+1}$ to $K_{n}$, {\it i.e.} they satisfy $f^{p_n}(K_n)=K_{n-1}$, $f^{a(I_n)}(\AA(I_n))=I_n$, and $f^{\sigma_\tau (K_{n})}(I_{n+1})=K_{n}$. 
Recall the statement of Lemma~\ref{l:doublement}~\,: 
\begin{lemi}[\ref{l:doublement}]
  $p_{n}\ge 2p_{n-1}$.
\end{lemi}
\proof Notice that $p_n=a(I_n)+\sigma_\tau (K_{n-1})$ and that  $p_{n+1}=a(I_{n+1})+\sigma(K_{n})$.

\begin{claimb}\label{c:4} The following inequality holds\,: $2r(I_n)\le a(I_n)\le (b+1)r(K_n)$.
\end{claimb}
\proof
The left  inequality comes from the definition of $\AA(I_n)$. Recall that  $\AA(I_n)$ is the pullback  of $W$ by $f^{b(I_n)}:\BB(I_n)\to I_n$, where $W$ is    $\LL_{f^{b(I_n)}(c_0)}(I_n)$.  Since  $\BB(I_n)\subset I_n$,   $b(I_n)\ge r(I_n)$.
Let $k$ be the time such that  $f^k(W)=I_n$, then   $k\ge r(I_n)$ since $W\subset I_n$. Therefore $ a(I_n)=k+b(I_n)\ge 2r(I_n)$.

The right inequality  is a corollary of Remark~\ref{r:returntime}  and Proposition~\ref{p:pteAB}. Indeed,  the point $c_0$
is contained in at most 
$b+1$ iterates of $\AA(I_n)$  in $\{f^i(\AA(I_n))\mid 0\le i\le a(I_n) \}$. Denote these iterates $f^{k_l}(\AA(I_n))$ with $k_0=0<k_1<\cdots<k_j=a(I_n)$. By Remark~\ref{r:returntime}.2), $k_{i+1}-k_i\le r(f^{k_i}(K_n))$, and since $f^{k_i}(K_n)\supset K_n$,  $ r(f^{k_i}(K_n))\le r(K_n)$, so that  $a(I_n)=\sum (k_{i+1}-k_i)\le (b+1)r(k_n)$ .

\begin{claimb}\label{c:5}  The following inequality holds\,: $ (2^{\tau+1}-2)r(K_n)\le \sigma_\tau(K_{n})\le 2r(I_{n+1})$.
\end{claimb}
\proof
 We prove first the right inequality. Let us  apply  Remark~\ref{r:returntime}2) at  two successors $A=\DD^{j+1}(K_n)$ and $A'=\DD^{j}(K_n)$\,:
we get $\sigma(\DD^{j}(K_n))\le r(\DD^{j+1}(K_n))$. Now successive applications of  Corollary~\ref{c:2successors}  give $\displaystyle r(\DD^{j+1}(K_n))\le \frac 1{2^{\tau-j-1}}r(\DD^\tau(K_n))$. Therefore,    $\displaystyle\sigma(\DD^{j}(K_n))\le \frac 1{2^{\tau-j-1}}r(I_{n+1})$. Hence $\displaystyle\sigma_\tau(K_n)=\sum_{j=0}^{\tau-1} \sigma(\DD^{j}(K_n)) \le \sum_{j=0}^{\tau-1}\frac 1{2^{\tau-j-1}}r(I_{n+1})$. Thus, $\sigma_\tau(K_n)\le 2r(I_{n+1})$. 

For the left inequality, note first that $\sigma(J)\ge 2r(J)$ by Corollary~\ref{c:2successors}. Applying this inequality successively, we get that $\displaystyle\sigma_\tau(K_n)=\sum_{j=0}^{j=\tau-1}\sigma(\DD^j(K_n))\ge \sum_{j=0}^{j=\tau-1} 2r(\DD^j(K_n))$. Then Corollary~\ref{c:2successors} implies that for all $j\ge 0$,  $r(\DD^j(K_n))\ge 2^{j+1}r(K_n)$. Therefore, $\sigma_\tau(K_n) \ge \sum_{j=0}^{j=\tau-1} 2^{j+1}r(K_n)= (2^{\tau+1}-2)r(K_n)$.

\vskip 2em
\noindent  {\it Proof of Lemma~\ref{l:doublement}.} By  claim~\ref{c:4} and~\ref{c:5} above, $p_n=a(I_n)+\sigma_\tau(K_{n-1})\le (b+1)r(K_n)+2r(I_{n})$. Thus $p_n\le (b+3)r(K_n)$ since $K_n\subset I_{n}$. The minorations coming from these two Lemmas
 imply that  $p_{n+1}=a(I_{n+1})+\sigma(K_{n})\ge 2r(I_{n+1})+(2^{\tau+1}-2)r(K_{n})\ge (2^{\tau+2}-2)r(K_{n}) $ since $I_{n+1}\subset K_n$. Hence,  the Lemma follows from the fact that for $\tau\ge b+1$ we have $(2^{\tau+2}-2)\ge 2(b+3)$.\cqfd

\section{Proof of Theorem~\ref{th:description}}\label{s:KU}
Up to replacing $f$ by some iterate, we can assume that the bounded component $U$ is fixed by $f$. We can also assume that $U$ contains only one critical point in it (up to doing some surgery). We  will also assume  that $K_U=K$\,:  
to achieve this  we use the fact that the restriction of $f$ to the domain bounded by an equipotential  around $K_U$ is a polynomial-like mapping. The Douady-Hubbard Theorem provides us with  a polynomial whose dynamics on the filled Julia set is conjugated to that of $f$ on $K_U$. We will still denote $f$ the new polynomial to which these three changes have possibly been performed.

By Theorem~\ref{th:fatouboundaryJordan}, the boundary of $\partial U$ is a Jordan curve. Hence, the Riemann map $\Phi_U: \D\to U$ extends continuously to the boundary 
as an homeomorphism  $\gamma_U:\S^1\to \partial U$ that  conjugates the dynamics (of the model map on $\S^1$ to $f$ on $\partial U$).
We 
 consider one of the graphs used in the article for constructing the puzzle. %
\begin{lem}\label{l:limb}Recall   the notation $End(z):=\Cap_{n\in \N}\ol{P_n(z)}$.
\begin{enumerate}
\item  If $z\in \partial U$ is eventually periodic then  either $End(z)=\{z\}$ or there exist two external rays $R_\infty(\zeta),R_\infty(\zeta')$ landing at $z$
and separating  $End(z)\setminus \{z\}$ from $\ol U$\,;
\item
 if $z\in \partial U$ is not eventually periodic  then $End(z)=\{z\}$\,;
 \item
 in both cases when $End(z)=\{z\}$, there exists at least one external ray converging to $z$. 
 \end{enumerate}
\end{lem}
\proof  When $z$ is periodic, this is Proposition~\ref{p:kiwi}. The case when $z$ is eventually periodic follows by taking pre-images. This proves 1).
When $z\in \partial U$ is not eventually periodic and satisfies $(\star)$; this is Corollary~\ref{c:star}. Lemma~\ref{l:cassimples} together with Corollary~\ref{c:nonpersist} and Corollary~\ref{c:star} imply that the only points $z\in \partial U$ where the local connectivity could fail are those such that $perCrit(x)$ is not empty. Those points were considered in details in section~\ref{s:persrec}, where Corollary~\ref{c:Ex} showed $End(x)=\{x\}$.
This proves 2). 

 For 3),  let $R_\infty(\zeta_n),R_\infty(\zeta'_n)$ be the two  external rays of $\partial P_n(z)$ landing on $\partial U$ and let $\zeta$ be   the limit of the sequence $(\zeta_n)$.  Then    the end of  $R_\infty(\zeta)$  enters every  puzzle $ P_n(z)$, so it converges to  $z=\cap_n P_n(z)$.
\cqfd
\begin{definition}For $t\in \S^1$,  let $\RR$ denote the union of all the external rays landing at $z=\gamma_U(t)$. Let $\widetilde U$ denote the component of $\C \setminus( \RR\cup\{z\})$ containing $U$. Then take $L_t:=K_U\cap (\C\setminus \widetilde U)$, it is called {\it a limb}.
\end{definition}

\begin{proprietes} The limbs have the following  properties\,:
\begin{itemize}
\item The set $L_t$ is connected\,: it follows from the fact that  $\gamma_U(t)$ belongs to  the closure of every component of $\C\setminus \RR$\,;
\item  $L_t\cap \ol U=\{\gamma_U(t)\}$, this follows from Lemma~\ref{l:limb}\,;
\item  $K_U= \Cup_{t\in \S^1} L_t\cup U $ by definition\,;
\item $L_t$ reduces to one point if and only if there is exactly one ray converging to $z$.
\end{itemize}
\end{proprietes}
\begin{lem}Let $z\in partial U$ be a point  such that    two rays  at least converge to $z=\gamma_U(t)$. Then either $L_t$ contains a critical point or  $L_t$ is mapped to some $L_{t'}$ which contains a critical point.
\end{lem}
\proof By definition, there exist two external rays $R_\infty(\zeta)$ and  $R_\infty(\zeta')$ which land at $z$ and  which separates $L_t\setminus \{z\}$ from $\ol U$. Assume that  $L_t$ does not contain a critical point. The sector $\C\setminus \widetilde U$ is a disk   between these two external rays, it contains no critical point (since the critical points are all in $K_U$), the boundary  is mapped homeomorphically  to   $R_\infty(D\zeta) \cup R_\infty(D\zeta') \cup \{f(z)\}$  (if the polynomial is of degree $D$). Thus it is mapped homeomorphically to the sector, which does not  contain $U$, and is bounded by  $R_\infty(D\zeta)$ and  $R_\infty(D\zeta')$. Therefore $f(L_t)=L_{dt}$ if the degree in $U$ is $d$. But the multiplication by $D$ will eventually cover $\S^1$. So some Limb of an iterated image of $z$ has to contain a critical point.\cqfd
\begin{lem}\label{l:8.4}
If $L_t$ contains a critical point then at least two external rays converge to $\gamma_U(t)$. \end{lem}
\proof If  the  critical point is $z=\gamma_U(t)$, then the pull back of any external ray landing at $f(z)$ consists in  two external rays landing at $z$. If the critical point is in $L_t\setminus \{z\}$,   then in particular $L_t$ is not reduced to a point. The result follows then from the definition.
\cqfd

Note that Corollary~1 follows from Lemma~\ref{l:8.4} since at  an eventually periodic point $x\in \partial U$ $End(x)=\{x\}$. So by construction $(P_n(x)\cap J)$ form a sequence of connected neighborhoods  of $x$ in $J$ whose diameter tends to $0$.

\end{document}